\documentclass[journal,twoside,web]{ieeecolor}
\usepackage{makecell}
\usepackage{bbm, dsfont}
\usepackage{comment}
\usepackage{stmaryrd}
\usepackage{generic}
\usepackage{cite}
\usepackage{amsmath,amssymb,amsfonts}
\usepackage{graphicx}
\usepackage{textcomp}
\usepackage{graphicx}
\usepackage{caption}
\usepackage{algorithmicx,algpseudocode,algorithm}
\usepackage{mathtools}

\usepackage{subcaption}

\newtheorem{theorem}{Theorem}
\newtheorem{lemma}{Lemma}

\newtheorem{assumption}{Assumption}

\newtheorem{remark}{Remark}

\def\BibTeX{{\rm B\kern-.05em{\sc i\kern-.025em b}\kern-.08em
    T\kern-.1667em\lower.7ex\hbox{E}\kern-.125emX}}

\newcommand{\col}{\text{col}}
\newcommand{\iset}[1]{\llbracket #1 \rrbracket}

\begin{document}
\title{Social optimum of a Linear Quadratic\\ Collective Choice Model under Congestion}
\author{Noureddine Toumi, Roland Malham\'{e}, 
and Jérôme Le Ny
\thanks{This work was supported by NSERC grant RGPIN-5287-2018, FRQNT grant 291809, and a fellowship from Polytechnique Montreal.}
\thanks{The authors are with the department of Electrical Engineering, 
Polytechnique Montreal and GERAD, QC H3T1J4, Montreal, Canada. Email:  
{\tt noureddine.toumi, roland.malhame,jerome.le-ny@polymtl.ca}}}
\maketitle

\begin{abstract}

This paper investigates the social optimum for a dynamic linear quadratic collective choice problem 
where a group of agents choose among multiple alternatives or destinations. The agents' common objective 
is to minimize the average cost of the entire population. A naive approach to finding a social optimum 
for this problem involves solving a number of linear quadratic regulator (LQR) problems that increases exponentially with the population size. 
By exploiting the problem's symmetries, we first show
that one can equivalently solve a number of LQR problems equal to the number of destinations,
followed by an optimal transport problem parameterized by the fraction of 
agents choosing each destination. 
Then, we further reduce the complexity of the solution search by defining 
an appropriate system of limiting equations, whose solution is used to obtain a strategy 
shown to be asymptotically optimal as the number of agents becomes large.
The model includes a congestion effect captured by a negative quadratic term 
in the social cost function, which may cause agents to escape to infinity in finite time. 
Hence, we identify sufficient conditions, independent of the population size, for the existence of the social optimum. Lastly, we investigate the behavior of the model through numerical simulations in different scenarios.
\end{abstract}


\section{Introduction}
Networked systems consist of highly interdependent dynamic components or agents that interact with each other
to achieve a given individual or collective purpose. These systems hold a significant importance in contemporary 
society, with applications spanning numerous domains including power distribution, transportation, epidemiological 
modeling, and analyzing social networks.
Social optimization refers to the analysis of networked systems where agents aim to optimize 
a collective outcome, i.e., find a social optimum. As the number of agents and the complexity 
of their interactions grow, identifying a social optimum becomes increasingly complex. 
To address this challenge, various techniques are employed. 

A first approach consists in using the model's properties, .e.g., its linear quadratic (LQ) structure
and symmetries, to derive an explicit analytical solution. 
Using this approach, Arabneydi et al. \cite{Arabneydi2016LinearSolutions} derive an 
explicit solution for a discrete-time LQ team optimization model under full and partial information structures, and Huang et al. \cite{Huang2021LinearControl} identify the solution of a social optimization problem with controlled diffusions.
Another frequently used method consists in using the person-by-person (PbP) optimality principle,
reframing the search for a social optimum as a non-cooperative game \cite{Yuksel2013StochasticSystems}.
Concretely, one solves the social optimization problem as a sequential optimization problem,
reducing it to a single-agent problem against a given population behavior.
In general, PbP optimality is only a necessary condition for social optimality.
Nevertheless, under specific conditions such as static games with convex costs 
and specific information structures, both solutions coincide 
\cite[Theorem 18.2]{vanSchuppenJanHandVilla2015CoordinationSystems}.

Another important framework for the approximate analysis of large interacting dynamic multi-agent systems
is the theory of Mean Field Games (MFG) \cite{Lasry2006b,Huang2007}.
It assumes that inter-agent couplings vanish asymptotically 
with the population size, so that, in the limit, unilateral deviation by any agent does not affect 
the global population behavior. The MFG and PbP approaches can be combined to address social optimization problems. This is done in \cite{Sen2016MeanCharacterization,Salhab2018,Huang2020MeanEquation}, 
where the authors derive suboptimal strategies converging to the social optimum as the population size increases. Nonetheless, this approach requires certain restrictive technical assumptions, as for PbP optimality.

This article focuses on a LQ social optimization model, with a particular focus 
on population transportation. For this application, modeling congestion is a critical factor to consider,
often resulting in a system of partial differential equations, which hinders the analytical
tractability of the solution \cite{Lasry2006b, Lachapelle2011, Aurell2018a}.

In this paper, we address congestion by introducing in the agents' cost functions a negative quadratic term that decreases with the agents' relative spacing, thus preserving the model's 
LQ structure. 
An advantage of LQ models is their computational tractability, which makes them suitable for real-time control and optimization. Nonetheless, the change 
in cost structure invalidates the conditions for PbP optimality used in \cite{Salhab2018}. 
Instead, we start from a finite population and exploit the model's symmetries and LQ structure 
as in \cite{Huang2021LinearControl} to derive the social optimum, as well as a more efficiently computable suboptimal solution approximating the social optimum 
as the population size increases.

Therefore, we consider a group of agents with continuous linear dynamics, initially spread in the Euclidean space $\mathbb{R}^n$. By the time horizon $T$, individual agents 
must have moved toward one of the alternatives within a set of destination points. This is to be achieved while minimizing an overall joint quadratic cost. 
When we fix the fraction of agents going to each destination, we show how to compute the optimal controls and individual agent choices by solving a number of LQR problems equal to 
the number of destinations followed an optimal transport (OT) problem parameterized by these fractions of agents.
Finally, considering the limiting system of optimality equations as the number of 
agents grows, we establish the convexity of the cost with respect to the vector 
of destination-bound agent fractions. This leads to an efficient optimization method
for the infinite population, as well as a policy for the finite agent system
that is shown to be asymptotically optimal as the number of agents becomes large.
 
The rest of the paper is organized as follows. 
In Section \ref{section:Problem Formulation}, we present the general formulation of our model. 
In Section \ref{section: solving the finite population problem}, we first establish 
the agents' optimal trajectories under an arbitrarily fixed agent-to-destination choices rule. Then, we provide an algorithm to identify a solution for our proposed model. 
In Section \ref{section:Riccati Escape Time}, we provide an upper bound on the time horizon $T$ that guarantees the  existence of a solution,
independently of the population size. 
In Sections \ref{section: limiting system of equations} and \ref{section: good aproximation}, we provide a limiting system 
of equations whose solution can be used to approximate the social optimum for 
a large population. In Section \ref{section: suboptimal algorithm}, we provide an algorithm to compute the latter solution. Section \ref{section: simulation} illustrates the results in
simulations, and Section \ref{section: conclusion} concludes the paper. 
In Table \ref{tab: notations}, we collect some of the notations used throughout the paper.

\begin{table}[h]
\caption{Notation table}
\centering
\begin{tabular}{|c|c|}
\hline
$n$& Dimension of an agent's state vector \\
\hline
$N$& Population size\\
\hline
$D$ &  Number of destinations \\
\hline
$\otimes$ & Kronecker product \\
\hline
$|.|$ & Set cardinality operator \\
\hline
$\llbracket m \rrbracket$ & The set $\{1,\ldots,m\}$, for $m$ a positive integer \\
\hline
$\Gamma^{N}$ &  A vector or matrix that depends on the population size $N$\\
\hline
$\Gamma^{(N)}$ &  \thead{Block diagonal matrix with $N$ repeated \\ blocks $\Gamma$, i.e., $I_{N} \otimes \Gamma$}\\
\hline
$1_{N}$& All-ones column vector of size $N$ \\
\hline
$1_{N N}$ & All-ones matrix of size $N \times N$ \\
\hline
$0_{n}$& All-zeros column vector of size $n$ \\
\hline
$0_{nn}$& All-zeros matrix of size $n \times n$ \\
\hline
$|y|_2$ & \thead{The Euclidean norm $\sqrt{y^Ty}$ for a vector $y$} \\
\hline
$[y]_R^2$ & \thead{The quadratic product  $y^TRy$, for a \\ given symmetric matrix $R$ and vector $y$}\\
\hline
$\text{col}(.)$ & Operator stacking its vector arguments into a column vector \\ 
\hline
$\Omega^N$ & $\begin{aligned}
&=\Bigg \{ P^N=\col \left( P_1^N,\hdots,P_{D}^N \right) \in \mathbb{R}^{D} \Big | \exists N_1,\hdots, N_{D} \in \mathbb{N}\\
& \sum_{j=1}^{D} N_j = N , \text{and } P_j^N = \frac{N_j}{N}, \text { for } 
j \in \iset{D} \Bigg \} 
\end{aligned}$ \\
\hline
$\Omega$ & $\begin{aligned}
= \Big \{ & P=\col(P_1,\hdots,P_{D}) \in \mathbb{R}^{D} \Big | \sum_{j=1}^{D} P_j = 1 \\
&\text { and } P_j \geq 0, \text { for } j \in \iset{D} \Big \} \nonumber \end{aligned}$ \\
\hline
\end{tabular}
\label{tab: notations}
\end{table}

\section{Problem Formulation }
\label{section:Problem Formulation}

We consider a dynamic cooperative game involving $N \geq 2$ agents, initially 
randomly distributed identically and independently  
over a subset $E$ of the Euclidean space $\mathbb{R}^n$
according to a distribution $\mathcal{P}_0$, whose support is included in $E$.
Before a given time horizon $T$, these agents should reach the neighborhood of one of $D$ possible destinations, denoted $d_i \in \mathbb R^n$, for $1 \leq i \leq D$. More precisely, the agents cooperate to minimize a common social cost $J^N_{\text{soc}}$, i.e., to solve \begin{align} 
& \qquad \min_{\lambda^N \in \iset{D}^N}
\inf\limits_{U(.)} J^N_{\text{soc}}(\lambda^N,U) := 
\frac{1}{N}\sum\limits_{i=1}^N J_{i,\lambda_i}(u_i,\bar{x}^N), 
\label{eq:social cost} \\
& s.t. \quad \frac{d}{dt} x_i(t) = A  x_i +B u_i(t), \;\; x_i(0)=x_i^0, \;\; 
i \in \iset{N},
\nonumber 
\end{align}
with $x_i(t) \in \mathbb{R}^{n}$ 
the state vector of agent $i$, $u_i(t) \in \mathbb{R}^{m}$ its control vector, $\lambda_i \in \iset{D}$
the index of its final destination choice,
$\lambda^N = (\lambda_1,\ldots,\lambda_N)$ the destination choices for the $N$ agents, 
and $\bar{x}^N(t)=\frac{1}{N}\sum_{i=1}^{N}x_i(t)$ the empirical population mean state.
We denote $U(t)=\text{col}(u_1(t),\hdots,u_N(t)) \in \mathbb{R}^{Nm}$.
The individual costs in \eqref{eq:social cost} are
\begin{align}
&
J_{i,\lambda_i}(u_i,\bar{x}^N) 
=  \int_{0}^{T} \frac{1}{2} \Big ( -[x_i- \bar{x}^N]_{R_x}^2+[x_i - d_{\lambda_i}]^2_{R_d} 
\nonumber \\ 
 &+[u_i]^2_{R_u} \Big ) d t + \frac{1}{2}[x_i(T)- d_{\lambda_i}]^2_{M}. \label{eq:individual  cost}
\end{align} 
The matrices $R_u$ and $M$ are symmetric positive definite, while the matrices $R_d,R_x$ are 
symmetric positive semidefinite. 
The first term in the running cost in \eqref{eq:individual  cost} is negative and captures a congestion 
avoidance effect. Indeed, since we have
\[
\sum_{i=1}^{N}[x_i- \bar{x}^N]_{R_x}^2=\frac{1}{2N}\sum_{i=1}^{N}\sum_{j=1}^{N}[x_{i}-x_{j}]_{R_x}^2,
\]
avoiding the population mean trajectory is equivalent to maximizing 
squared interdistances between agents. 
The second term in the running cost increases whenever the distance between the $i^{th}$ 
agent position $x_i$ and its chosen destination $d_{\lambda_i}$ increases, representing a stress effect felt by individuals as long as they remain  far from their 
destination. Overall,
we aim to identify the optimal destination choices $\lambda^N$ and control strategies 
$U$ for the agents to collectively minimize \eqref{eq:social cost}. Note that since the running cost 
may be negative, agents may escape to infinity in finite time. 
Therefore, in our analysis, we provide in Section \ref{section:Riccati Escape Time} a sufficient condition for the existence of an optimal solution over the interval $[0,T]$. 

\section{Solving the Problem\\ for a Population of Finite Size}
\label{section: solving the finite population problem}

In this section, we provide a solution to the problem formulated in Section \ref{section:Problem Formulation}. 
The problem \eqref{eq:social cost} involves a minimization over 
both the destination choices $\lambda^N$ and the vector of control strategies $U$. 
Our approach consists in first addressing the optimization over $U$, while fixing $\lambda^N$. 
Then, we fix only the fraction of agents choosing each of the destinations $P^N$
and find the corresponding vector of optimal individual destination
choices $\lambda^N$. 
This step-wise approach enables us to identify the optimal social cost for all possible 
fixed vectors of agent fractions over destinations. At that stage, the identification of the optimal $\lambda^N$ 
and $U$ becomes simply a parameter optimization problem over all feasible vectors of
agent fractions over destinations.

\subsection{Agent Control Strategies with Given Agent-to-Destination Choices} 
\label{section:Agents Optimal trajectories}

We assume for now that the agents' final destinations $\lambda^N$ are fixed and known. 

We let $d_{\lambda^N}=\text{col}(d_{\lambda_1},\hdots,d_{\lambda_N})$, and
\begin{equation} \label{eq: def of Q^N}
Q^N \coloneqq I_N \otimes (R_d-R_x) + \frac{1}{N} 1_{N N} \otimes R_x.
\end{equation}
The following lemma, proved in Appendix \ref{proof: social cost manipulation}, 
expresses the social cost \eqref{eq:social cost} as a function of the global state
$X=\text{col}(x_1,\hdots,x_N) \in \mathbb{R}^{Nn} $, which satisfies
\[
\frac{d}{dt} X= A^{(N)} X +B^{(N)} U.
\]

\begin{lemma}
\label{lemma: social cost augmented state}
The social cost \eqref{eq:social cost} can be expressed as follows
\begin{align}
    &J^N_{\text{soc}}(\lambda^N,U) =        \label{eq:social cost augmented state} \\
    &\int_0^T \frac{1}{2N} \Big \{
    [X]^2_{Q^N} - 2 d_{\lambda^N}^T R_d^{(N)} X  
    +[d_{\lambda^N}]^2_{R_d^{(N)}}+[U]^2_{R_u^{(N)}}  \Big \} dt \nonumber \\
    &+ \frac{1}{2N} [X(T)-d_{\lambda^N}]^2_{M^{(N)}}. \nonumber 
\end{align} 
\end{lemma}
\vspace{0.2cm}

To find a control strategy $U^*$ minimizing
\eqref{eq:social cost augmented state}, 
we consider the quadratic value function
\begin{equation} \label{eq: value function}
V^N(t,\lambda^N) = \frac{1}{2N}X^T(t)\varphi^N(t) X(t)+\frac{1}{N}\Psi^N(t)^T X(t)+\chi^N(t),
\end{equation}
with $\varphi^N(t) \in \mathbb{R}^{Nn \times Nn}$, $\chi^N(t)\in \mathbb{R}$,
and
\begin{equation}
\Psi^N(t)=\col(\Psi^N_1(t),\hdots,\Psi^N_N(t)) \in \mathbb{R}^{Nn}, \label{eq:big psi}
\end{equation}
where $\Psi^N_i(t) \in  \mathbb{R}^{n}$ for all $i \in \iset{N}$.

By applying the dynamic programming, we obtain the following lemma, which we prove in Appendix \ref{proof: dynamic programming}
\begin{lemma}
\label{lemma: dynamic programming}
The time-varying coefficients $\varphi^N(t)$, $\Psi^N(t)$, and $\chi^N(t)$  should satisfy the following 
differential equations
\begin{align}
    \frac{d}{dt}\varphi^N &= (\varphi^N)^T S^{(N)} \varphi^N-\varphi^N A^{(N)}-(A^{(N)})^T\varphi^N -Q^N,    
    \label{eq:riccati all population} \\
    \frac{d}{dt}\Psi^N &= \left(\varphi^N S^{(N)}-{A^{(N)}}^T \right)\Psi^N+R_d^{(N)}d_{\lambda^N}, 
    \label{eq:psi all population} \\
    \frac{d}{d t}\chi^N &= \frac{1}{2N} \left( [\Psi^N]^2_{S^{(N)}}-[d_{\lambda^N}]^2_{R_d^{(N)}} \right), 
    \label{eq:chi bar} 
\end{align} 
with the terminal conditions
\begin{align}
\varphi^N(T) &= M^{(N)}, \; \Psi^N(T)=-M^{(N)}d_{\lambda^N},  \label{eq:phi Psi terminal}  \\
\chi^N(T) &= \frac{1}{2N}[d_{\lambda^N}]^2_{M^{(N)}}, \label{eq:chi bar terminal} 
\end{align}
and
\begin{align}
    S=B R_{u}^{-1} B^{T}.
\end{align}
Moreover, the optimal control \( U^*(t) \) is given by
\begin{align}
    U^*(t)=-(R_u^{(N)})^{-1}(B^{(N)})^T \left( \varphi^N(t) X(t) +\Psi^N(t) \right). 
    \label{eq:control all population}
\end{align}
\end{lemma}
We note that $\Psi^N$ and $\chi^N$ are both solutions over $[0,T]$ of linear differential equations, 
which depend on $\varphi^N$.
Hence, if a  solution exists for the backward Riccati equation 
\eqref{eq:riccati all population} over [0,T], 
then, both $\Psi^N$ and $\chi^N$ are well defined. However, the matrix $Q^N$ is not necessarily positive semidefinite. 
Therefore, the Riccati equation \eqref{eq:riccati all population} may exhibit a finite escape time. 
First, from \cite[Theorem 5.3]{Rami2002}, we have the following lemma.

\smallskip
\begin{lemma}   \label{lemma: finite optimal cost}
The optimal cost $\inf\limits_{U(.)} J^N_{\text{soc}}(\lambda^N,U)$ is finite
for any $\lambda^N \in \iset{D}^N$ if and only if
the Riccati equation \eqref{eq:riccati all population} has a solution defined over $t \in [0,T]$.
\end{lemma}
\smallskip

In Section \ref{section:Riccati Escape Time}, we shall identify sufficient conditions, 
independent of the population size $N$, that guarantee the existence of a finite solution 
for \eqref{eq:riccati all population} over [0,T]. 
For the rest of this section, we assume that such conditions are met. 
Then, we have the following result,
proved in Appendix \ref{proof: varphi decomposition}. 

\smallskip
\begin{lemma}
\label{lemma: varphi decomposition}
The solution $\varphi^N$ to the Riccati equation \eqref{eq:riccati all population} can be decomposed as 
\begin{equation}
\label{eq:varphi decomposition}
\varphi^N(t)=I_N \otimes \varphi^N_1(t) + \frac{1}{N} (1_{N N}-I_N)\otimes \varphi^N_2(t),
\end{equation}
where the diagonal and off-diagonal components $\varphi^N_1(t) \in \mathbb R^{n \times n}$ 
and $\varphi^N_2(t) \in \mathbb R^{n \times n}$ satisfy the equations
\begin{align}
    \frac{d}{dt}\varphi^N_1 &= [\varphi^N_1]^2_S+\frac{N-1}{N^2}[\varphi^N_2]^2_S-\varphi^N_1 A-A^T\varphi^N_1  \nonumber \\ 
    & \;\; - R_d + \frac{N-1}{N}R_x, \quad \varphi^N_1(T)=M,  \label{eq:diagonal component for the riccati} \\
    \frac{d}{dt}\varphi^N_2 &= (\varphi^N_1)^T S \varphi^N_2+(\varphi^N_2)^TS\varphi^N_1+\frac{N-2}{N}[\varphi^N_2]^2_S \nonumber \\
    & \;\; -\varphi^N_2 A -A^T\varphi^N_2-R_x, \quad \varphi^N_2(T)=0_{nn}. \label{eq:off-diagonal component for the riccati}
\end{align}
\end{lemma}
\medskip

For any given $j \in \iset{D}$,
we denote by $N_j$ the number of agents 
whose final target is $d_j$ and by $P_j^N \coloneqq \frac{N_j}{N}$
the fraction of agents selecting the destination $j$.
Note that $P_{D}^N=1-\sum\limits_{j=1}^{D-1} P_j^N$ necessarily.
We also define $P^N=\col(P_1^N,\hdots,P_{D}^N) \in \Omega^N$,
with $\Omega^N$ defined in Table \ref{tab: notations}.
In the next lemma, proved in Appendix \ref{proof: psi decomposition}, we show that 
if two agents choose the same destination, then their corresponding components 
of the solution $\Psi^N$ of \eqref{eq:psi all population} are identical.\\
\begin{lemma} \label{lemma: psi decomposition} 
Define $\psi^N_j(t,P^N) \in \mathbb R^{n}$, $j \in \iset{D}$, as the solutions to the following system of equations
\begin{align}
\frac{d}{d t} &\psi^N_{j}(t,P^N) = \left( \varphi^N_1S-A^T-\frac{1}{N}\varphi^N_2S \right) \psi^N_{j}+R_d d_j 
\label{eq:psi principle components} \\
& \;\; +\varphi^N_2 S \left( \sum\limits_{k=1}^{D} P_k^N \psi^N_{k} \right), \quad \psi^N_{j}(T,P^N)=-Md_j. \nonumber
\end{align} 
Then, recalling (\ref{eq:big psi}), we have $\Psi^N_i \equiv \psi^N_{\lambda_i}$ for $i \in \iset{N}$.
\end{lemma}
\smallskip
Next, we identify an analytical form for the solution of \eqref{eq:psi principle components}.
The proof of this result is given in Appendix \ref{proof: psi form finite}.
\smallskip
\begin{lemma}
\label{lemma: psi form finite}
For any $j \in \iset{D}$, the function  $\psi^N_{j}(t,P^N)$ is
affine with respect to $P^N$ and satisfies the following relation 
\begin{equation}
\label{eq: linear psi}
\psi^N_{j}(t,P^N)= \left( \sum_{k=1}^{D-1} P^N_k \alpha^N_k(t) \right) -\beta^N_j(t),
\end{equation}
where $\alpha^N_k(t) \in \mathbb R^{n}$ and $\beta^N_j(t) \in \mathbb R^{n}$, 
are solutions of\\
\begin{align}
    \frac{d}{dt} \alpha^N_k &= \left( \varphi^N_1 S - A^T + \frac{N-1}{N} \varphi^N_2 S \right) \alpha^N_k 
    -\varphi^N_2 S (\beta^N_k-\beta^N_D) \nonumber \\
    & \;\; \alpha_k^N(T)=0_{n}, \;\; \forall k \in \iset{D-1}, 
    \label{eq: finite alpha differential equation} \\
    \frac{d}{dt} \beta^N_j &= \left( \varphi^N_1 S - A^T -\frac{1}{N}\varphi^N_2 S \right) \beta^N_j + 
    \varphi^N_2 S \beta^N_D - R_d d_j, 
    \nonumber \\
    & \;\; \beta^N_j(T)=M d_j, \;\; \forall j \in \iset{D}. 
    \label{eq: finite beta differential equation}
\end{align} 
\end{lemma}
\smallskip
In the following, we define
\begin{align} 
\begin{split}   \label{eq: alphaN matrix}
& \alpha^N=\begin{bmatrix} \alpha_1^N & \ldots & \alpha_{D-1}^N & 0_n \end{bmatrix} \in \mathbb R^{n \times D}, \\
& \beta^N=\begin{bmatrix} \beta_1^N &\ldots &\beta_D^N \end{bmatrix} \in \mathbb R^{n \times D}. 
\end{split}
\end{align}
Finally, in the next lemma, proved in Appendix \ref{proof: chi form finite}, 
we establish the form of the solution of \eqref{eq:chi bar}. 
\begin{lemma}
\label{lemma: chi form finite}
The solution $\chi^N(t,P^N)$ of \eqref{eq:chi bar} is a quadratic form of $P^N$ and we have that

\begin{align}
& \frac{d}{dt} \chi^N =(P^N)^T W^N P^N +\frac{1}{2}\sum_{j=1}^{D}P_j^N\left([\beta_j^N]^2_{S}-[d_j]^2_{R_d}\right),     \nonumber   \\
&\chi^N(T,P^N)=\frac{1}{2} \sum_{j=1}^{D} P_j^N[d_j]^2_M, \label{eq: chi polynomial finite} 
\end{align} 
where $W^N(t) \in \mathbb R^{D \times D}$ is defined as
\begin{align}   \label{eq: W N}
W^N(t) = \left( \frac{1}{2}\alpha^N(t)-\beta^N(t) \right)^T S \, \alpha^N(t). \end{align} 
\end{lemma}
\smallskip
Using Lemmas \ref{lemma: varphi decomposition} and \ref{lemma: psi form finite}, as proved in Appendix \ref{proof:value function rewritten}, we can now rewrite 
the value function $V^N(t,\lambda^N)$ in \eqref{eq: value function} as follows.
\begin{lemma}
The value function in \eqref{eq: value function} can be expressed as
\label{lemma:value function rewritten}
\begin{align}
&V^N(t,\lambda^N) =\frac{1}{2N}\sum\limits_{i=1}^{N}x_i^T(t)\left(\varphi^N_1(t)-\frac{\varphi^N_2(t)}{N}\right)x_i(t) \nonumber \\
&+\frac{1}{2}(\bar{x}^N(t))^T\varphi^N_2(t)\bar{x}^N(t)+\frac{1}{2N}\sum\limits_{i=1}^{N}|x_{i}-\beta_{\lambda_i}^N|_2^2 \nonumber \\
&+\sum\limits_{j=1}^{D-1}P^N_j(\alpha_j^N(t))^T\bar{x}^N(t)-\frac{1}{2}\sum\limits_{j=1}^{D}P^N_j|\beta_{j}^N(t)|_2^2 \nonumber \\
& -\frac{1}{2N}\sum\limits_{i=1}^{N}|x_i(t)|_2^2+\chi^N(t,P^N). \label{eq:simplified value function}
\end{align}
\end{lemma}
\smallskip
Denote by $u_{\lambda_i}^{N,*}$ the optimal control of an arbitrary agent with index $i \in \llbracket N \rrbracket$,
destination $\lambda_i$, and state vector $x_i \in \mathbb R^n$. 
Using \eqref{eq:control all population} and Lemma \ref{lemma: varphi decomposition}, the agent's optimal control is
\begin{align}  
u_{\lambda_i}^{N,*}(t,x_i,P^N) &=-(R_u)^{-1}B^T \left( \left( \varphi^N_1-\frac{\varphi_2^N}{N} \right) x_i \right. \nonumber\\
&\left.\vphantom{\int_1^2}+\varphi^N_2\bar{x}^N(t,P^N)+\psi^N_{\lambda_i}(t,P^N) \right) ,\label{eq:social optimal control}
\end{align}
with
\begin{align}
\frac{\partial }{\partial t} \bar{x}^N(t,P^N)&=\left[A-S\left(\varphi^N_1+\frac{N-1}{N}\varphi^N_2\right)\right]\bar{x}^N(t,P^N) \nonumber \\
&-S \sum_{j=1}^D P_j^N \psi^N_{j}(t,P^N) \label{eq: mean dynamics}.
\end{align}

\subsection{Complete Problem Solution}
\label{section: finite assignment problem} 

To identify the optimal solutions $U^*$ and $\lambda^N$ to problem \eqref{eq:social cost} for 
a free agent-to-destination choice, we follow a two-step procedure. First, we identify 
the optimal agent-to-destination choice for any given final destination occupancy probability 
vector. Subsequently, we use the result of this first step to turn the search 
of the optimal solution $U^*$, $\lambda^N$ into a parameter optimization problem over the set $\Omega^N$.

Assume first that we fix the destination occupation probability vector $P^N \in \Omega^N$. 
When evaluating \eqref{eq:simplified value function} at $t=0$, all terms
are then fixed
except for $\sum_{i=1}^{N} |x_{i}(0)-\beta_{\lambda_i}^N(0)|^2$. Therefore, identifying the agent-to-destination 
choices $\lambda^N$, for a given probability vector $P^N \in \Omega^N$ reduces to solving 
the following optimal transport (OT) problem \cite{Peyre2019ComputationalTransport}
\begin{align}
\label{eq: discrete Monge prime}
& C^N(P^N)=\min_{\lambda^N \in \iset{D}^N} \frac{1}{N}\sum_{i=1}^{N} |x_{i}(0)-\beta_{\lambda_i}^N(0)|_2^2 \\
& \text{s.t. } \frac{| \{ i \in \iset{N} | \lambda_i=j  \}|}{N}=P_j^N, \; \forall j \in \iset{D}.  \nonumber
\end{align} 
Using the integrality theorem \cite[Theorem 14.2]{Vanderbei2020LinearExtensions}, 
Problem \eqref{eq: discrete Monge prime} can be solved exactly and efficiently through 
its linear programming relaxation, which we introduce also later in \eqref{eq: discrete Kantorovich prime}.

For every probability vector $P^N \in \Omega^N$, once we solve \eqref{eq: discrete Monge prime},
we can define the associated optimal social cost function
\begin{align}
\label{eq:simplified value function 2}
&J^N(P^N) = \frac{1}{2N} \sum_{i=1}^{N}x_i^T(0)\left(\varphi^N_1(0)-\frac{\varphi^N_2(0)}{N}\right)x_i(0) \\
&+\frac{1}{2} \left( \bar{x}^N(0) \right)^T \varphi^N_2(0)\bar{x}^N(0)+\frac{1}{2}C^N(P^N) \nonumber \\
&+\sum\limits_{j=1}^{D-1}P^N_j(\alpha_j^N(0))^T\bar{x}^N(0)-\frac{1}{2}\sum\limits_{j=1}^{D}P^N_j|\beta_{j}^N(0)|_2^2
\nonumber \\
& -\frac{1}{2N}\sum\limits_{i=1}^{N}|x_i(0)|_2^2+\chi^N(0,P^N), \nonumber 
\end{align}
which is yet to be optimized over $P^N \in \Omega^N$.
To identify the optimal controls $U^*$ and vector $\lambda^N$ for problem \eqref{eq:social cost}, 
we can then use Algorithm \ref{algo:finite population}. 
However, Algorithm \ref{algo:finite population} is computationally expensive 
for many agents and destinations because it searches exhaustively for the optimal $P^N$ 
in the set $\Omega^N$, which is of size $\left(\begin{array}{c}N+D-1 \\ D-1 \end{array}\right)$.
In Section \ref{section: limiting system of equations}, we propose a related limiting 
system of equations leading to a cost convex in the probability vector $P^N$, 
which becomes arbitrarily close to $J^N(P^N)$ as the population size $N$ goes to infinity.
This ultimately leads to a more efficient search procedure for an asymptotically optimal
final destination distribution.

\begin{algorithm}[htbp]
\caption{Computing the optimal $U^*$ and $\lambda^N$ for \eqref{eq:social cost}.}
\label{algo:finite population}
\begin{algorithmic}[1]
\State Compute solutions $\varphi_1^N$, $\varphi_2^N$, $\alpha_k^N, k \in \iset{D-1}$, 
and $\beta_j^N, j \in \iset{D}$ of \eqref{eq:diagonal component for the riccati}, 
\eqref{eq:off-diagonal component for the riccati}, \eqref{eq: finite alpha differential equation}, 
\eqref{eq: finite beta differential equation}.
\For{each $\underline P^N \in \Omega^N$}
\State Compute $\chi^N$ using \eqref{eq: chi polynomial finite}
\State Compute $C^N(\underline P^N)$ in \eqref{eq: discrete Monge prime} using linear programming
\State Compute the cost $J^N(\underline P^N)$ in \eqref{eq:simplified value function 2}
\EndFor
\State Identify the vector $P_{opt}^{N} \in \Omega^N$ minimizing $J^N$
\State Retrieve $C^N(P_{opt}^N)$ and $\lambda^N$ solving \eqref{eq: discrete Monge prime}. 
This $\lambda^N$ corresponds to the optimal agent-to-destination choices.
\State Compute the agents' optimal controls in \eqref{eq:social optimal control} using the optimal agent-to-destination choices $\lambda^N$.
\end{algorithmic}
\end{algorithm}

\subsection{Existence of a Solution to the Riccati Equation}
\label{section:Riccati Escape Time}

The results of Sections \ref{section:Agents Optimal trajectories} and \ref{section: finite assignment problem} 
assume that a solution exists for the Riccati equation \eqref{eq:riccati all population}
over the whole interval $[0,T]$ when integrated backwards. 
However, this backward differential equation may exhibit a finite escape time 
$\Delta_{esc}$
such that the solution ceases to exist for $t \leq T - \Delta_{esc}$
\cite{abou2012matrix,TOUMI2024111420}, and moreover it is possible that
$T - \Delta_{esc} > 0$.
Note also that $\Delta_{esc}$ depends on $Q^N$ 
defined in \eqref{eq: def of Q^N}, and hence on $N$.
Our goal in this section is to identify a sufficient condition,
independent of $N$, for the existence of a solution over [0,T].

For any given matrices $C,D,F,G$ and time horizon $T$, we denote by 
$\textit{HDRE}(T,C,D,F,G)$ the Hermitian differential Riccati equation of the form
\begin{equation} \begin{aligned}
\frac{d}{dt} \varphi = \varphi^T C \varphi - \varphi D - D^T \varphi - F, \;\; \varphi(T)=G.
\nonumber \end{aligned}\end{equation} 
Hence, equation \eqref{eq:riccati all population} is $\textit{HDRE}(T,S^{(N)},A^{(N)},Q^N,M^{(N)})$. 
We also denote by $\Delta_{esc}(C,D,F,G)$ the (possibly infinite) escape time of the 
equation, i.e., the maximal $\Delta > 0$ such that the solution exists on $(T-\Delta,T]$.
The following theorem is proved in Appendix \ref{proof:bounding escape time}.

\begin{theorem} \label{thm: bound on time horizon}
For any $T < \Delta_{esc}(S,A,R_d-R_x,M)$, the Riccati equation \eqref{eq:riccati all population} has a solution 
over the time interval $[0,T]$.
\end{theorem}

For the rest of the paper, we make the following assumption, which guarantees by Lemma \ref{lemma: finite optimal cost}
that the optimal cost \eqref{eq:social cost} is finite for all $N \geq 1$.
\begin{assumption}
The time horizon T is strictly smaller than $\Delta_{esc}(S,A,R_d-R_x,M)$, which is independent of $N$.
\end{assumption}

Finally, to identify  
$\Delta_{esc}(S,A,R_d-R_x,M)$, we first define the operator $\Delta(t,\varphi_0)$, 
where $\varphi_0 \in \mathbb{C}^{n \times n}$, by
\begin{equation} \begin{aligned}
\label{eq:escape time criterion}
\Delta(t,\varphi_0) &= \det \left[I+\int_{0}^{t} e^{\breve{A} p} S e^{\breve{A}' p} d p \; (M-\varphi_0)\right]\\
\breve{A}(\varphi_0)&=A-S\varphi_0, \breve{A}'(\varphi_0)=A^T-\varphi_0S. 
\end{aligned}\end{equation} 
Assume the Riccati equation $\textit{HDRE}(T,S,A,R_d-R_x,M)$ 
admits an equilibrium $\varphi_0 \in \mathbb{C}^{n \times n}$. For instance, it is the case when $(A,B)$ is controllable \cite[Theorem 16]{Lancaster1980ExistenceEquation}.
Using the results in \cite{Sasagawa1982OnEquations}, we establish the following lemma, 
which also provides a method to compute $\Delta_{esc}(S,A,R_d-R_x,M)$.
\begin{lemma} \label{lemma: general escape criterion} 
Let $\varphi_0$ be an equilibrium of the Riccati equation $\textit{HDRE}(T,S,A,R_d-R_x,M)$.
The escape time $\Delta_{esc}(S,A,R_d-R_x,M)$
corresponds to the first time $t>0$ where $t \rightarrow \Delta_d(t,\varphi_0)$ vanishes. 
Moreover, this time does not depend on the choice of equilibrium $\varphi_0$.
\end{lemma}
For the specific case where $A=0$ and the matrices $S,R_x,R_d$ and $ M $ are diagonal, 
an explicit analytical expression of $\Delta_{esc}(S,A,R_d-R_x,M)$ is identified 
in \cite[Lemma 8]{TOUMI2024111420}.

\section{The Limiting System of Equations}
\label{section: limiting system of equations}

In this section, we define a limiting system of equations by letting $N$ 
go to infinity in the equations of Sections \ref{section:Agents Optimal trajectories} 
and \ref{section: finite assignment problem}. We associate with this system an adequately 
defined cost function $J(P)$ where vector $P \in \Omega$ plays in the limit $N \to \infty$ 
the role of $P^N \in \Omega^N$. We prove in Section \ref{sec: continuum strategy} that the agent-to-destination choices 
resulting from minimizing $J(P)$ as well as the proposed limit control strategy
can serve as good approximation for the optimal strategy of Section \ref{section: solving the finite population problem}. 
This leads to an $\varepsilon-$optimal control strategy for Problem \eqref{eq:social cost}, 
with $\varepsilon$ going to zero as $N$ goes to infinity. We call this control strategy the \emph{continuum strategy}.

An advantage of the continuum strategy is that each agent can compute its control input by only knowing the initial population distribution $\mathcal{P}_0$, and not the exact initial positions of all agents
or even the number $N$ of agents. This reduces information exchanges and can lead to a more practical solution.
Additionally, \eqref{eq: discrete Monge prime} is replaced by a semi-discrete OT problem \cite{Peyre2019ComputationalTransport},
which is independent of the number $N$ of agents.
Finally, continuous optimization methods are available to optimize $J$.
In fact, we show in Section \ref{section: asymptotic solution properties},
under some assumptions on the problem parameters,
that the cost $J$ is convex, which leads us to propose in Section \ref{section: suboptimal algorithm} a new algorithm
that is significantly more efficient than the exhaustive search over $\Omega^N$ in 
Algorithm \ref{algo:finite population}, for a small loss of optimality.
In Section \ref{sec:limiting system}, we define the limiting system of equations. 
Then, given a probability vector $P \in \Omega$, we solve in Section \ref{sec: limiting Agent-to-destination Assignement} 
a limiting OT problem between $\mathcal P_0$ and $P$, whose optimal value is used in the definition of the cost 
function $J(P)$. In Section \ref{section: asymptotic solution properties}, we establish the continuity and convexity of $J$.
Before proceeding, we impose from now on the following condition, in order to guarantee that
all quantities introduced are well defined and simplify the discussion of the OT problem.
\begin{assumption}  \label{asspt: E compact}
The set $E$ is compact.
\end{assumption}
In particular, a consequence of Assumption \ref{asspt: E compact} is that all moments of $\mathcal P_0$ are finite.

\subsection{Defining the Limiting System of Equations}
\label{sec:limiting system}

To obtain the limiting system of equations, 
starting from \eqref{eq:diagonal component for the riccati}, 
\eqref{eq:off-diagonal component for the riccati}, define $\varphi_1$, $\varphi_2$ solutions of
\begin{align}
\frac{d}{dt}\varphi_1 &= [\varphi_1]^2_S - \varphi_1 A - A^T \varphi_1 - (R_d - R_x), 
\label{eq:diagonal component for the riccati inf} \\
\varphi_1(T) &= M, \nonumber \\
\frac{d}{dt}\varphi_2 &= \varphi_1^T S \varphi_2 + \varphi_2^T S \varphi_1 + [\varphi_2]^2_S - \varphi_2 A - A^T \varphi_2 - R_x, \nonumber \\
\varphi_2(T) &= 0_{nn}.     \label{eq:off-diagonal component for the riccati inf}
\end{align}

Next, based on \eqref{eq: finite alpha differential equation}, 
\eqref{eq: finite beta differential equation}, \eqref{eq: alphaN matrix}, define
\begin{align*}     
\alpha(t) &=\begin{bmatrix} \alpha_1(t)& \ldots & \alpha_{D-1}(t) & 0_n\end{bmatrix} \in \mathrm{R}^{n \times D},   \\
\beta(t) &=\begin{bmatrix} \beta_1(t)& \ldots &\beta_D(t) \end{bmatrix} \in \mathrm{R}^{n \times D},
\end{align*}
with
\begin{align}
\frac{d}{dt} \alpha_j &= 
(\varphi_1S-A^T+\varphi_2S) \alpha_j-\varphi_2S(\beta_j-\beta_D), 
\label{eq: alpha differential equation} \\
\alpha_j(T)&=0_{n}, \quad \forall j \in \iset{D-1}, \nonumber \\
\frac{d}{dt} \beta_j &= (\varphi_1S-A^T) \beta_j +\varphi_2S\beta_D - R_d d_j, 
\label{eq: beta differential equation} \\
\beta_j(T) &= Md_j, \quad \forall j \in \iset{D}.  \nonumber
\end{align} 
Finally, from \eqref{eq: W N}, define for all $t \geq 0$
\begin{equation} 
W(t) = \left(\frac{1}{2}\alpha(t)-\beta(t)\right)^T S \alpha(t) \label{eq: W inf}.
\end{equation}

Based on \eqref{eq: chi polynomial finite}, define for all 
$P =col(P_1,\hdots,P_{D}) \in \Omega$, the solution $\chi(\cdot,P)$ of
\begin{align}
\frac{d}{dt}\chi &= 
P^T W P +\frac{1}{2}\sum_{j=1}^{D}P_j\left([\beta_j]^2_{S}-[d_j]^2_{R_d}\right), \label{eq: chi polynomial infinite} \\
\chi(T,P) &= \frac{1}{2} \sum_{j=1}^{D} P_j [d_j]^2_M, \nonumber
\end{align}
Then, based on \eqref{eq: linear psi}, \eqref{eq:social optimal control} and \eqref{eq: mean dynamics}, 
the control for an agent at state $x \in \mathbb R^n$ with final target $d_{j}$, for $j \in \{1,\hdots,D \}$, is
\begin{align}  \label{eq:approx control}
&\tilde{u}_{j}(t,x,P) =-(R_u)^{-1}B^T \left(  \varphi_1 x 
+ \varphi_2\bar{x}(t,P)  +\psi_{j}(t,P) \right),
\end{align}
with
\begin{align}
\label{eq: inf psi}
\psi_{j}(t,P)&= \left( \sum_{k=1}^{D-1} P_k \alpha_k(t) \right) -\beta_j(t),\\
\frac{\partial }{\partial t} \bar{x}(t,P)&=\left[A-S\left(\varphi_1+\varphi_2\right)\right]\bar{x}(t,P)-S \sum_{j=1}^D P_j \psi_{j}, \nonumber\\
\bar{x}(0,P) &= \int_{x \in E} x \, d \mathcal P_0. \label{eq: x_bar}
\end{align} 
Finally, consider instead of \eqref{eq: discrete Monge prime} the following OT problem
\begin{equation} \begin{aligned}
\label{eq:semi-discrete kantorovich prime}
C(P)=
&\inf \limits_{\gamma \in \mathcal{M}(E\times \iset{D})} \sum_{j=1}^D \; \int\limits_{x \in E} |x-\beta_{j}(0)|_2^2 \; d \gamma(x,j) \\
& s.t. \qquad \gamma(E,j) = P_{j}, \; \forall j \in \iset{D}, \\
& \qquad \sum_{j=1}^{D}\gamma(A,j)=\mathcal{P}_0(A), \; \forall A \in \mathcal{B}(E),
\end{aligned} \end{equation} 
where $\mathcal{M}(E \times \iset{D})$ is the set of joint probability measures 
on $E \times \iset{D}$ and $\mathcal{B}(E)$ are the $\mathcal P_0$-measurable subsets of $E$. 

Corresponding to \eqref{eq:simplified value function 2}, we define the following cost function 
\begin{align}
&J(P) = \frac{1}{2}\int\limits_{x \in E}x^T\varphi_1(0)x \, d \mathcal{P}_0
+\frac{1}{2}\bar{x}(0)^T\varphi_2(0)\bar{x}(0)+\frac{1}{2}C(P) \nonumber \\
&+\sum\limits_{j=1}^{D-1}P_j\alpha_j(0)^T\bar{x}(0)-\frac{1}{2}\sum\limits_{j=1}^{D}P_j|\beta_{j}(0)|_2^2 
-\frac{1}{2}\int\limits_{x \in E}|x|_2^2 \, d \mathcal{P}_0 \nonumber \\
&+\chi(0,P), 
\label{eq:simplified value function infinite}
\end{align}
with the definition $\bar{x}(0) = \int_{x \in E} x \, d \mathcal P_0$ 
for the expected initial condition.
Note that all integrals in \eqref{eq:simplified value function infinite} are well defined under Assumption \ref{asspt: E compact}.

To compute the value of $J(P)$ for a given vector $P \in \Omega$, 
we need to solve the OT problem \eqref{eq:semi-discrete kantorovich prime},
as described in the next section. This solution yields a mapping that we use in 
Section \ref{section: asymptotic solution properties} to approximate the 
optimal agent-to-destination choices when the population size $N$ is large.

\subsection{Solving The OT Problem for a Given Vector $P$}
\label{sec: limiting Agent-to-destination Assignement}

Here we recall some features of the (semi-discrete) OT problem \eqref{eq:semi-discrete kantorovich prime},
see, e.g., \cite[Chapter 5]{Peyre2019ComputationalTransport} for additional details and references.
For a fixed $P \in \Omega$, this is an OT problem with quadratic cost between the measure $\mathcal P_0$ 
and the discrete measure $\sum_{j=1}^N P_j \delta_{\beta_j(0)}$, where $\delta_x$ denotes the Dirac measure
at $x \in \mathbb R^n$, i.e., for any set $A \subset \mathbb R^n$, $\delta_x(A) = 1$ if $x \in A$ and $\delta_x(A) = 0$ otherwise.
Even though \eqref{eq:semi-discrete kantorovich prime} is an infinite-dimensional linear
program, it is known that the infimum is attained and, by Kantorovitch duality \cite{Villani2009OptimalNew}, 
its value is equal to the supremum of the dual function $J^D(\cdot,P): \mathbb R^D \to \mathbb R$, defined as 
\begin{align}
&J^D(g,P) \coloneqq 
\int_{ E } \min\limits_{j \in \iset{D} }(|x-\beta_j(0)|_2^2-g_j) \, d \mathcal{P}_0(x) \label{eq:dual} \\
& \quad \quad \quad \quad \quad + \sum\limits_{j=1}^{D} P_{j} g_{j}, 
\;\; \text{ with } g \coloneqq \col(g_{1},\hdots,g_{D}). \nonumber
\end{align} 
To describe the optimal solution of \eqref{eq:semi-discrete kantorovich prime} and simplify the discussion, 
it is convenient to impose the following additional condition on $\mathcal P_0$, which is satisfied 
for instance when $\mathcal P_0$ is absolutely continuous with respect to the Lebesgue measure.
\begin{assumption}  \label{asspt: measure zero}
Hyperplanes in $\mathbb R^n$ have $\mathcal P_0$-measure zero.
\end{assumption}

For any given vector $g \in \mathbb{R}^D$, the term $\min\limits_{j \in \iset{D}}(|x-\beta_j(0)|_2^2-g_j)$ 
in the dual function \eqref{eq:dual} induces a partition of the set $E$ into polygonal cells 
$\mathcal{C}_{j}(g), j \in \iset{D}$, defined as

\begin{align}   \label{eq: general cells}
&\mathcal{C}_{j}(g) \coloneqq \\
&\Big\{ x \in E: |x-\beta_j(0)|_2^2-g_j \leq |x-\beta_k(0)|_2^2-g_k, \forall k \in \iset{D} \Big \}. \nonumber
\end{align}
This partition is called a power diagram, a type of generalized Voronoi diagram, 
see, e.g., \cite{Aurenhammer1998Minkowski-TypeClustering} for more details.
Under Assumption \ref{asspt: measure zero}, the dual function can also be rewritten
\begin{align*}
J^D(g,P) &= \sum_{j=1}^D \int_{\mathcal C_j(g)} |x-\beta_j(0)|_2^2 \; d \mathcal P_0(x) \\
& \quad + \sum_{j=1}^D g_j (P_j - \mathcal P_0(\mathcal C_j(g)).
\end{align*}
The following lemma is proved in Appendix \ref{proof: dual function lemma}.
\begin{lemma}   \label{lemma: dual optimum}
Let $P \in \Omega$. The function $g \mapsto J^D(g,P)$ defined in \eqref{eq:dual} is concave.
Moreover, under Assumption \ref{asspt: E compact}, it admits a maximizer $g^* \in G$ with $$G=\{ g \in \mathbb R^D | \; |g|_2\leq \sup\limits_{\substack{x \in E \\ j \in \llbracket D \rrbracket }}|x-\beta_{j}(0)|_2^2 \}.$$
\end{lemma}
Given $P \in \Omega$, take $g^*(P) \in \arg \max_g J^D(g,P)$, whose existence is guaranteed by Lemma \ref{lemma: dual optimum}.
We then define the corresponding polygonal cells
\begin{equation} \label{eq:cells}
\mathcal{C}_{j}^*(P):=\mathcal{C}_{j}(g^*(P)), \;\; j \in \iset{D},
\end{equation} 
and the mapping $\lambda: E \times \Omega \rightarrow \iset{D}$, such that  
\begin{equation} \label{eq: mapping}
\lambda(x,P) = \min \{ j \in \iset{D}: x \in \mathcal{C}^*_j(P)\}.
\end{equation}
The mapping $\lambda(\cdot,P)$ sends each $x$ in cell $C^*_j(P)$ to destination $j$,
and resolves ambiguities at the cell boundaries by arbitrarily assigning the cell with the smallest index.
Under Assumption \ref{asspt: measure zero}, the rule used to choose a destination for the agents
initially on the cell boundaries has no impact on the cost $C(P)$, since these initial conditions
occur with probability zero.
Note also that $\lambda(\cdot,P)$ induces a (deterministic) plan $\gamma \in \mathcal M(E \times \iset{D})$ 
for the OT problem \eqref{eq:semi-discrete kantorovich prime} 
(more precisely, the measure $\gamma_\lambda = (\text{id},\lambda)_{\#} \mathcal P_0$, 
pushforward of $\mathcal P_0$ by $(\text{id},\lambda(\cdot,P)): E \mapsto E \times \iset{D}$).

We then have the following result \cite{Peyre2019ComputationalTransport}.
\begin{lemma}   \label{lemma: SDOT result}
Under Assumptions \ref{asspt: E compact} and \ref{asspt: measure zero},
the mapping $\lambda(\cdot,P)$ defined in \eqref{eq: mapping} induces an optimal solution for the OT problem
\eqref{eq:semi-discrete kantorovich prime}.
\end{lemma}
In Section \ref{sec: continuum strategy} we use $\lambda$ to approximate 
the optimal agent-to-destination choices $\lambda^N$.

\begin{remark}  \label{remark: finite dual optimum}
If we replace the compactness condition of Assumption \ref{asspt: E compact} by the weaker assumption that $\mathcal P_0$
has a finite second moment, the cost function $J$ can still be defined by \eqref{eq:simplified value function infinite}
and the OT problem \eqref{eq:semi-discrete kantorovich prime} still admits a minimizer plan $\gamma^*$, 
however the description of this minimizer becomes more cumbersome, in particular to handle situations where 
some destinations $j \in \iset{D}$ have $P_j = 0$. In that case, it is possible
that an optimal dual vector $g^*$ maximizing \eqref{eq:dual} necessarily has some components $g_j$ equal 
to $-\infty$. This can be seen for example by taking $n=1$, $D=2$, $\mathcal P_0$ to be a Gaussian distribution
on $\mathbb R$, and $P = [1,0]^T$, in which case $\mathcal C_2^*(P)$ must be empty and if $g_1$ is finite, 
this can only be achieved by taking $g_2 = -\infty$. The discussion could be generalized to take such issues
into account by removing unassigned destinations, however this also complicates the presentation of 
the algorithm of Section \ref{section: suboptimal algorithm}, as $P$ can be on the boundary of $\Omega$. 
Hence, Assumption \ref{asspt: E compact} allows us to simplify the presentation by guaranteeing the existence
of a finite dual optimum in $\mathbb R^D$, for any $P \in \Omega$.
\end{remark}

\subsection{Continuity and Convexity of $J$}
\label{section: asymptotic solution properties}
We now establish the continuity and convexity of the function $J$ defined
in \eqref{eq:simplified value function infinite}. 
For this,
it is sufficient to verify the continuity and convexity of the functions 
$P \rightarrow C(P)$ and $P \rightarrow \chi(0,P)$, since the other terms of $J$ are either linear in $P$
or do not depend on $P$. 
The next lemma is proved in Appendix \ref{proof: convexity of C(P)}.
\begin{lemma}
\label{prop: psi cost part is convex}
Under Assumption \ref{asspt: E compact}, the function $P \rightarrow C(P)$ is continuous and convex.
\end{lemma}

Furthermore, \eqref{eq: chi polynomial infinite} shows
that $P \rightarrow \chi(0,P)$ is a quadratic form in $P$, thus continuous. 
As a result, $J$ is continuous.
To prove the convexity of $\chi(T,P)$, we consider the following simplifying assumption.
\begin{assumption}
\label{assmp: diagonal matrices}
The matrices $A,B,R_x,R_d,R_u,M$ are diagonal matrices.
\end{assumption}
Then, under this assumption, we have the following lemma, which we prove 
in Appendix \ref{proof: convexity of chi}.
\begin{lemma}
\label{prop: chi is convex} 
Under Assumption \ref{assmp: diagonal matrices}, the function $P \rightarrow \chi(0,P)$ is convex
on $\Omega$.
\end{lemma}

Combining the results in Lemmas \ref{prop: psi cost part is convex} and \ref{prop: chi is convex}, 
we obtain the following theorem.
\begin{theorem} \label{thm: convexity of J}
Under Assumptions \ref{asspt: E compact} and \ref{assmp: diagonal matrices}, the cost $J$ defined in 
\eqref{eq:simplified value function infinite} is a continuous and convex function on $\Omega$.
\end{theorem}

\section{The Continuum Strategy \\as a Suboptimal Strategy}
\label{sec: continuum strategy}

As explained in the introduction of Section \ref{section: limiting system of equations},
under the \textit{continuum strategy}, the agents use the limiting control \eqref{eq:approx control} 
as well as the agent-to-destination choices \eqref{eq: mapping} associated 
with a probability vector $P$ minimizing the cost function $J$.
In the following, we first show
in Section \ref{section: good aproximation} that the continuum strategy leads 
to an $\varepsilon-$optimal solution for Problem \eqref{eq:social cost}, with $\varepsilon$ 
going to zero as $N$ goes to infinity. Then, in Section \ref{section: suboptimal algorithm}, 
we propose a second algorithm (Algorithm \ref{algo:inf population}) to compute the $\varepsilon$-optimal strategy,
which is computationally more efficient than
Algorithm  \ref{algo:finite population}.
Since we rely on the limiting OT problem to define the continuum strategy, 
we work under Assumptions \ref{asspt: E compact} and \ref{asspt: measure zero}.

\subsection{Near-Optimality of the Continuum Strategy}
\label{section: good aproximation}

In this section, we show that the social cost incurred when using 
the continuum strategy approaches the optimal cost.
Let 
\begin{equation}
 J^N_{opt}=\min\limits_{P \in \Omega^N}J^N(P) \text{ and }  \label{eq: opt cost def}
J_{opt}=\min\limits_{P \in \Omega}J(P),
\end{equation}
with $J^N$ and $J$ defined in \eqref{eq:simplified value function 2} and 
\eqref{eq:simplified value function infinite} respectively.
Our proof begins by establishing that $J^N_{opt}$ converges almost surely to $J_{opt}$
as $N \to \infty$.
Then, we show that the social cost incurred when employing the continuum strategy 
also converges almost surely to $J_{opt}$, which leads us to the desired conclusion.

First, define $\bar{J}^N$ as the extension of $J^N$ over $\Omega$, i.e.,
\begin{equation} \begin{aligned}    \label{eq: J bar N}
&\bar{J}^N(P)=\frac{1}{2N}\sum\limits_{i=1}^{N}x_i^T(0)\left(\varphi^N_1(0)-\frac{\varphi^N_2(0)}{N}\right)x_i(0) \\
&+\frac{1}{2}\bar{x}^N(0)^T\varphi^N_2(0)\bar{x}^N(0)+\frac{1}{2}\bar{C}^N(P) \\
&+\sum_{j=1}^{D-1}P_j(\alpha_j^N(0))^T\bar{x}^N(0)-\frac{1}{2}\sum\limits_{j=1}^{D}P_j|\beta_{j}^N(0)|_2^2\\
& -\frac{1}{2N}\sum\limits_{i=1}^{N}|x_i(0)|_2^2+\chi^N(0,P),
\end{aligned} \end{equation}
for any vector $P \in \Omega$, with
\begin{equation} \begin{aligned}
\label{eq: discrete Kantorovich prime}
&\bar{C}^N(P)=\min\limits_{\gamma \in \mathcal{M}( \iset{N} \times \iset{D} ) } 
\sum_{\substack{i \in \iset{N} \\ j \in \iset{D} }}  |x_i(0)-\beta^N_{j}(0)|_2^2 \gamma(i,j)\\
\end{aligned}\end{equation} 
\begin{equation*} \begin{aligned}
s.t. &\sum_{i=1}^{N}\gamma(i,j) = P_{j}, \forall j \in \iset{D}, 
\sum_{j=1}^{D}\gamma(i,j)=\frac{1}{N}, \forall i \in \iset{N},  \\
&\gamma(i,j) \geq 0, \forall (i,j) \in \iset{N} \times \iset{D},
\end{aligned}\end{equation*} 
where $\mathcal{M}(\iset{N} \times \iset{D})$ is the set of joint probability mass functions 
on $\iset{N} \times \iset{D}$ and $\chi^N(0,P)$ is obtained by replacing 
$P^N \in \Omega^N$ by the vector $P \in \Omega$ in \eqref{eq: chi polynomial finite}.
The following lemma is proved in Appendix \ref{proof: value function uniform convergence}. 
\begin{lemma}   \label{lemma: value function uniform convergence}
Under Assumption \ref{asspt: E compact}, with  probability one, the function $\bar{J}^N$ converges uniformly to ${J}$ on $\Omega$.
\end{lemma}
 This leads us to the next lemma, proved in Appendix \ref{proof: conv optimal cost to limit}, 
 which concludes the first step of our proof.
\begin{lemma}   \label{lemma: conv optimal cost to limit}
Under Assumption \ref{asspt: E compact},
with the definitions \eqref{eq: opt cost def}, the optimal social cost $J^N_{opt}$ converges almost surely to $J_{opt}$
as $N \to \infty$.
\end{lemma}
Recall that the continuum strategy is associated with a vector $P$ minimizing the cost function $J$,
and the destination choices \eqref{eq: mapping} for that specific $P$. When this strategy is applied 
for a finite population of size $N$, the actual proportion of agents choosing destination $j \in \iset{D}$ is
\begin{align}
F^N_j(P) \coloneqq \frac{1}{N} \Big| \{i \in \iset{N} \; | \; x_i(0) \in \mathcal{C}_j^*(P)\} \Big|, \label{eq: F^N_j}
\end{align}
which differs from $P_j$.
Define the operator $F^N: \Omega \rightarrow \Omega^N$ by 
\begin{align}
F^N(P) = \col \left( F^N_1(P),\hdots,F^N_D(P) \right). \label{eq: F^N}
\end{align}
Next, for a given vector $P \in \Omega$, if each agent $i \in \iset{N}$ uses the agent-to-destination 
choice $\lambda(x_i(0),P)$ given by \eqref{eq: mapping} and the control $\tilde{u}_{\lambda(x_i(0),P)}$
from \eqref{eq:approx control} leading to a state trajectory denoted $\tilde{x}_i(t,P)$,
they incur an individual cost $J_{i,\lambda(x_i(0),P)}$ defined in \eqref{eq:individual  cost}
and hence a social cost
\begin{align}
\tilde{J}^N(P) &= \frac{1}{N}\sum\limits_{i=1}^N J_{i,\lambda(x_i(0),P)}(\tilde{u}_{\lambda(x_i(0),P)},\tilde{x}^N), \label{eq: tilde J}
\end{align} 
where $\tilde{x}^N(t,P) = \frac{1}{N} \sum_{i=1}^N \tilde{x}_i(t,P)$ denotes the corresponding agents' mean trajectory, which satisfies the differential equation
\begin{align}
\frac{\partial }{\partial t}\tilde{x}^N(t,P)&=\left(A-S\varphi_1\right)   \tilde{x}^N(t,P)- S \varphi_2 \bar{x}(t,P) \nonumber\\
&-S\sum_{j=1}^{D}F^N_j(P)\psi_{j}(t,P),  \label{eq: x_tilde}\\
\text{with } &\tilde{x}^N(0,P)=\frac{1}{N}\sum_{i=1}^{N} x_i(0). \nonumber 
\end{align}
Recall the definition of $J_{opt}$ in \eqref{eq: opt cost def}, 
let $P^* \in \arg\min\limits_{P \in \Omega}J(P)$, 
and define $\tilde{J}^N_{subopt}=\tilde{J}^N(P^*)$.
We prove the following lemma in Appendix \ref{proof: conv sub-optimal cost to limit}.
\begin{lemma} \label{lemma: conv sub-optimal cost to limit}
Under Assumptions \ref{asspt: E compact} and \ref{asspt: measure zero},
as $N \to \infty$, the cost $\tilde{J}^N(P)$ converges almost surely to $J(P)$, 
for all $P \in \Omega$. Hence, $\tilde{J}^N_{subopt}$ converges almost surely to $J_{opt}$ as $N \to \infty$.
\end{lemma}

Since for all $N \in \mathbb{N},$ we have with probability one
\[ 
|J^N_{opt}-\tilde{J}^N_{subopt}| \leq |J^N_{opt}-J_{opt}|+|\tilde{J}^N_{subopt}-J_{opt}|,
\]
the following theorem follows from 
Lemmas \ref{lemma: conv optimal cost to limit} and \ref{lemma: conv sub-optimal cost to limit}.
\begin{theorem} \label{theorem: epsilon optimum}
Under Assumptions \ref{asspt: E compact} and \ref{asspt: measure zero},
the continuum strategy is an $\varepsilon-$optimal strategy for the $N$-agent population problem, 
with $\varepsilon$ converging almost surely to $0$ as $N$ goes to infinity.
\end{theorem}

\subsection{Numerical Computation of the Continuum Strategy}
\label{section: suboptimal algorithm}
In this section, we provide a new algorithm (Algorithm \ref{algo:inf population}) to compute 
the continuum strategy. The algorithm consists of two nested loops. The outer loop minimizes 
the cost $J$ defined by \eqref{eq:simplified value function infinite} over $\Omega$ 
using the projected subgradient method \cite{zaslavski2020projected}. 
  
For a given probability vector $P$, the inner loop uses the projected gradient ascent method to find a
maximizer $g^*(P)$ of $J^D$ and the corresponding cells $\mathcal{C}_{i}^*(P)$ defined in \eqref{eq:cells}. 
Therefore, we start 
by establishing the analytical expressions for all gradients required for the implementation of the proposed algorithm. 
We work throughout this section under Assumptions \ref{asspt: E compact} and \ref{asspt: measure zero}. Under the additional Assumption \ref{assmp: diagonal matrices}, by Theorem \ref{thm: convexity of J}
the cost $J$ is convex and hence using the result in \cite[Proposition 3.2.7]{bertsekas2015convex}, we get that the projected subgradient method with the dynamic step size rule \cite[Eq 3.21]{bertsekas2015convex} identifies a global minimizer 
$P$ of $J$ to define the continuum strategy, so that the approximation result of 
Theorem \ref{theorem: epsilon optimum} applies. First, the following lemma recalls a well-known analytical expression for the supergradient $\nabla_g J^D(g,P)$ of the 
concave dual function \eqref{eq:dual}, used by the inner loop of Algorithm \ref{algo:inf population}, 
see, e.g., \cite{Aurenhammer1998Minkowski-TypeClustering}.
A proof can be found in \cite[Appendix B]{JLN:TAC13:deployment}.
\begin{lemma}
Fix $P \in \Omega$ and let $g \in \mathbb R^D$.
The vector $\nabla_g J^D(g,P)$ with components
\begin{align}   \label{eq: dual gradient}
[\nabla_g J^D(g,P)]_j= P_j-\int\limits_{x \in \mathcal{C}_{j}(g)} d \mathcal{P}_0(x), 
\;\;\; \forall j \in \iset{D},
\end{align}
is a supergradient of $g \mapsto J^D(g,P)$ at $g$.
\end{lemma}
In the next lemma, proved in Appendix \ref{proof: subgradient}, we provide an expression for a subgradient of $J$.

\begin{lemma}
\label{Lemma: subgradient} 
Let $P \in \Omega$, and let $g^*(P) \in \mathbb R^D$ be a maximizer for \eqref{eq:dual}. 
Then the vector 
\begin{align}   \label{eq: gradient cost}
J'(P)= g^*(P) - W_{\text{int}} \, P+ H,
\end{align} 
is a subgradient of $J$ at $P$, where \begin{align}
W_{\text{int}}&= 2 \int_0^T W(t) \, dt, \label{eq: W_int}
\end{align}
and $H \in \mathbb R^D$ has components
\begin{align}  \label{eq: H}
H_j &= -\frac{1}{2}  \int_0^T [\beta_j]^2_S -[d_j]^2_{R_d} dt + \alpha_j(0)^T\bar{x}(0) \\
& \quad + \frac{1}{2} \left( [d_j]^2_{M}-|\beta_{j}(0)|_2^2 \right), \;\; j \in \iset{D}, \nonumber
\end{align}
with the convention $\alpha_D(0) = 0$.
\end{lemma} 

Algorithm \ref{algo:inf population} describes the steps for the computation of the continuum strategy. 
In this algorithm, we use the Euclidean projections over the probability simplex $\Omega$ and the set $G$ 
denoted by $Proj_\Omega(.)$ and $Proj_G(.)$ respectively. The implementation details for this operation can be found in \cite{wang2013projection}
and its complexity is $\mathcal{O}(D \log(D))$.
\begin{algorithm}[htbp]
\caption{Computing the continuum strategy for \eqref{eq:social cost}.}
\label{algo:inf population}
\begin{algorithmic}[1]
\Statex \textbf{Inputs:} initial $P^{out}$, $g$; stepsizes $s_{out}$, $s_{in}$; threshold $\delta$
\State Compute solutions $\varphi_1$, $\varphi_2$, $\alpha_i, i \in \iset{D-1}$, $\beta_i, i \in \iset{D}$, $W_{int}$, 
and $H$ from \eqref{eq:diagonal component for the riccati inf}, \eqref{eq:off-diagonal component for the riccati inf}, \eqref{eq: alpha differential equation}, \eqref{eq: beta differential equation}, \eqref{eq: W_int}, \eqref{eq: H}
\State $c_{out} \gets 0$
\Repeat
\Repeat
\State Compute $P^{in}_i \gets \mathcal{P}_0(C_j(g)), \forall j \in \iset{D}$, with $C_i(.)$ defined in \eqref{eq: general cells}
\State $g \gets g +  s_{in} (P^{out}-P^{in})$ \; (see \eqref{eq: dual gradient})
\Until{$|P^{in}-P^{out}|_{\infty} \leq \delta$}
\State $P^{out}_{old} \gets P^{out}$
\State Update the step size $s_{out}$ using the update rules of Equations 3.21, 3.22, and 3.23 in \cite{bertsekas2015convex}
\State $P^{out} \gets P^{out}-s_{out} (g - W_{\text{int}}P^{out}+ H)$ 
(see \eqref{eq: gradient cost}).
\State $P^{out} \gets Proj_\Omega(P^{out})$
\Until{$|P^{out}-P^{out}_{old}|_{\infty} \leq \delta$}

\State Using the final weight vector $g$, compute the partition $\mathcal{C}_{i}(g)$  
(see \eqref{eq: general cells})
and the corresponding agent-to-destination mapping $\lambda$ (see \eqref{eq: mapping}).
\State Compute the agents' approximate controls \eqref{eq:approx control} using the 
mapping $\lambda$ from \eqref{eq: mapping} and $P=P^{out}$.
\end{algorithmic}
\end{algorithm}

\section{Numerical Illustration of the Model}
\label{section: simulation}
\begin{figure}[t!]
    \centering
    \begin{subfigure}[b]{0.45\linewidth}
        \centering
        \includegraphics[width=\linewidth]{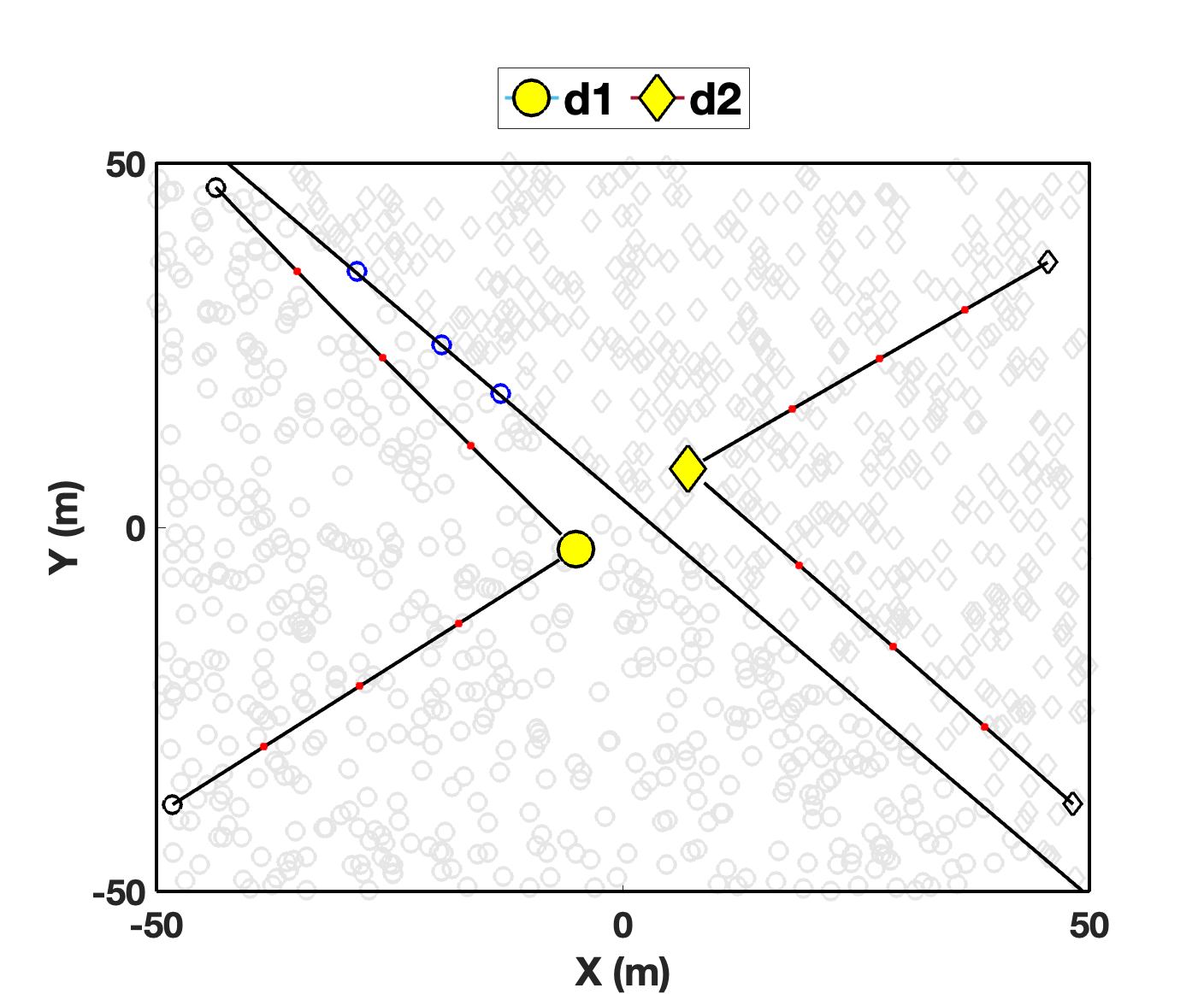}
        \caption{$T=3$.}
    \end{subfigure}
    ~ 
    \begin{subfigure}[b]{0.45\linewidth}
        \centering
        \includegraphics[width=\linewidth]{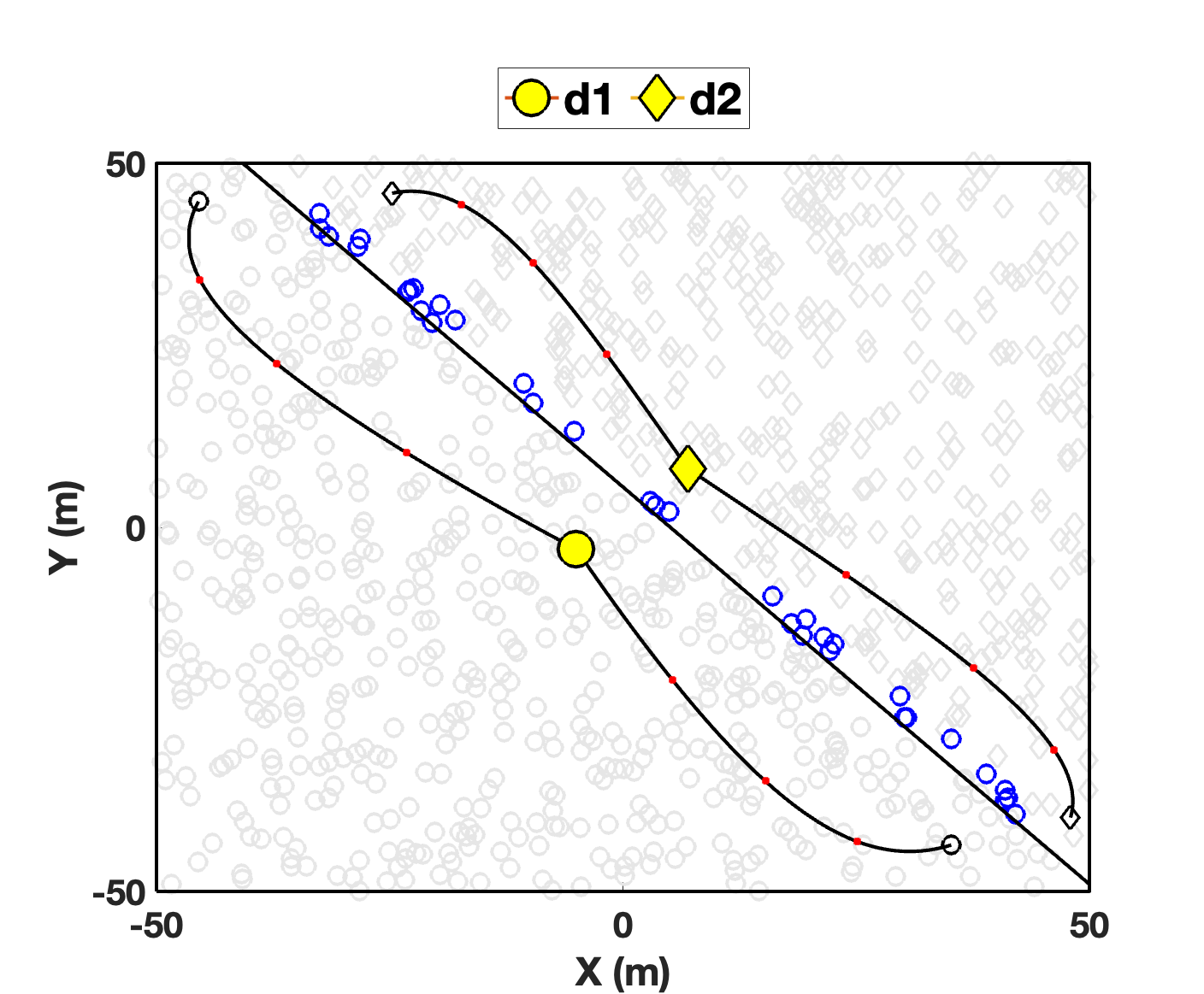}
        \caption{$T=10$.}
    \end{subfigure}
   \caption{Agent trajectories and destinations choices, for $1000$ agents. Destinations are depicted by the yellow, black-contoured shapes. The central black line delineates destination choices under the continuum strategy. Grey, blue, and black markers indicate the agents' initial positions, with their shapes corresponding to their chosen targets when employing the optimal strategy. 
The agents marked in blue are those whose optimal destination choice differs 
from that of the continuum strategy. Four agents and their corresponding trajectories are represented in black. Each trajectory displays three red dots, denoting the agents' positions at instants $T/4$, $T/2$, and $3T/4$, respectively.
}
\label{fig: trajectories}
\end{figure}

In this section, we compare the optimal strategy of Section \ref{section: solving the finite population problem} 
and the continuum strategy of Section \ref{sec: continuum strategy}. Hence, we simulate a group of agents initially spread over the domain $E = [-50,50] \times [-50,50]$ 
according to the uniform distribution. 
The agents move according to the dynamics $A= 0_{2 \times 2}$, $B=I_{2\times2}$. Two final destinations are located at the coordinates $d_1=[-5,-3]^T$ and $d_2=[7,8]^T$.
The stress effect matrix is $R_d = 0.1 \, I_{2 \times 2}$ and the social interaction effect matrix 
is $R_x = I_{2 \times 2}$. The other model parameters are  
$M= 400 \, I_{2\times2}$ and $R_u=50 \, I_{2\times2}$. Using \cite[Lemma 8]{TOUMI2024111420},
we obtain an upper bound on the time horizon $\Delta_{esc}(S,A,R_d-R_x,M)=23.4135$, so the Riccati equation
\eqref{eq:riccati all population} has a solution over the time interval $[0,T]$.
The input parameters for Algorithm \ref{algo:inf population} are $P^{out}=[0.5,0.5]^T$, $g=[0,0]^T$, 
$s_{in}=3500$, and $\delta=5.10^{-5}$. 
On Fig.\ref{fig: trajectories}, we simulate, for two different time horizons $T=3$ and $T=10$, 
one thousand agents and display their destination choices when deploying each strategy. 
Also, we plot the trajectories of four randomly chosen agents. Note that if an agent selects 
the same final destination for both strategies, its trajectories will be nearly identical 
for both cases. Therefore, we only represent one trajectory for each of the agents. 
 
 First, we note that agents trajectories are much more curved when $T= 10$ compared to those obtained for $T=3$, which can be seen as a precursor of the finite escape time phenomenon. A similar effect occurs when increasing the matrix $R_x-R_d$, or decreasing the matrix $R_u$, causing agents to increase their relative distances to avoid congestion. Also, we see that only the agents starting close to the boundary of the cells corresponding to the limiting OT problem make suboptimal destination choices.

\begin{figure}[t!]
    \centering
    \begin{subfigure}[b]{0.45\linewidth}
        \centering
        \includegraphics[width=\linewidth]{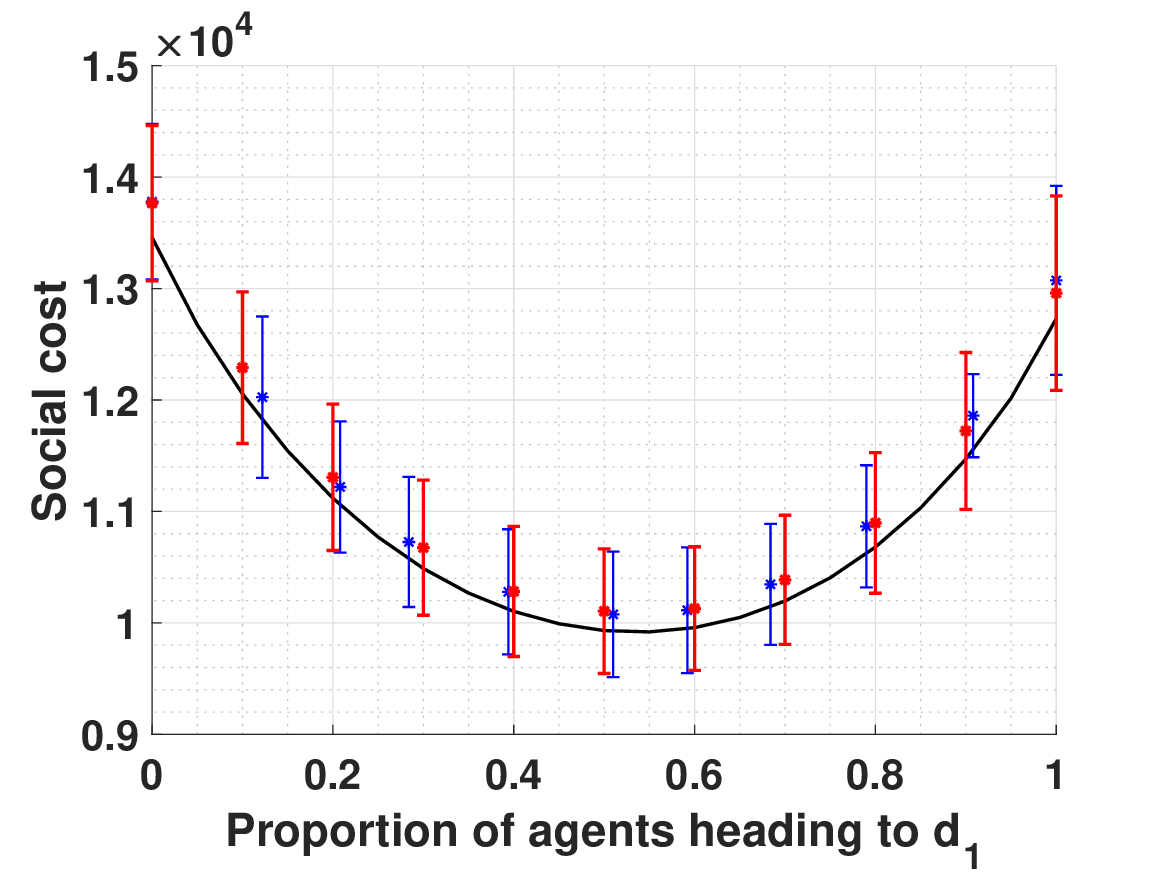}
        \caption{$100$ agents.}
    \end{subfigure}
    ~ 
    \begin{subfigure}[b]{0.45\linewidth}
        \centering
        \includegraphics[width=\linewidth]{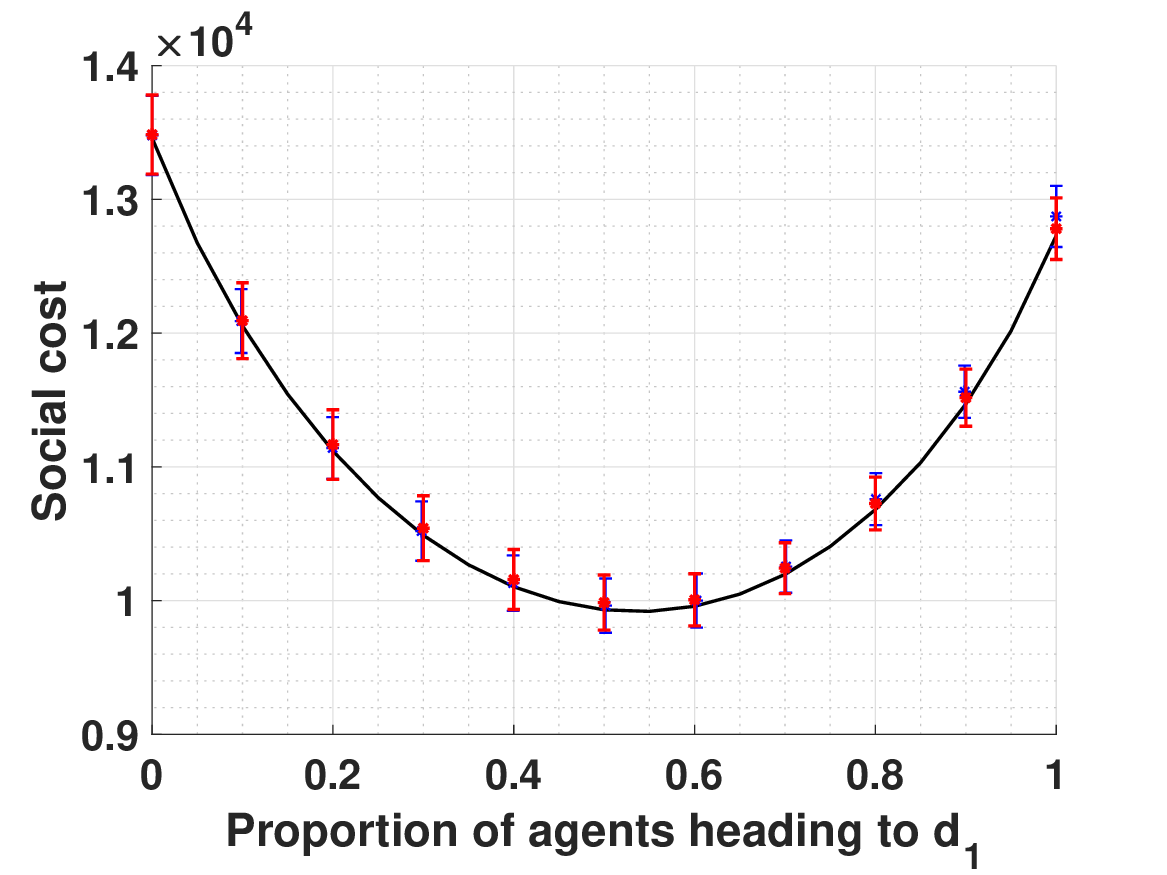}
        \caption{$1000$ agents.}
    \end{subfigure}
    \caption{Costs $J$ (\eqref{eq:simplified value function infinite}, black curve) 
    $\tilde{J}^N$ (\eqref{eq: tilde J}, in blue) and $J^N$ (\eqref{eq:simplified value function 2}, in red)
    for two populations.
    }
\label{fig: social costs}
\end{figure}

On Fig.~\ref{fig: social costs}, for two different populations sizes $100$ and $1000$, we plot the cost $J$ defined in \eqref{eq:simplified value function infinite}, the average and standard deviation of the costs $J^N$ and $\tilde{J}^N$ from \eqref{eq:simplified value function 2}, \eqref{eq: tilde J}, computed over 20 simulations and for a time horizon $T=3$. Running the Algorithm \ref{algo:inf population} gives that The optimal probability vector minimizing the cost $J$ is $[0.5535,0.4465]^T$. The costs $J^N$ and $\tilde{J}^N$ are evaluated on a probability grid with a discretization step of $0.1$. However, since agents employing the continuum strategy agents use the cells \eqref{eq:cells}, the resulting agent-to-destination proportion $F^N(P)$ in \eqref{eq: F^N} may exhibit slight deviations from the vector $P$. To reflect this, we choose to plot the average of the cost $\tilde{J}^N(P)$ against the average of the probability vectors $F^N(P)$, even though the cost $\tilde{J}^N(P)$ is computed as a function of $P$.

 First, note that the cost function $J$ is a convex function, supporting the statement of Theorem \ref{thm: convexity of J}. In addition, we observe that the standard deviation of the costs averages of $J^N$ and $\tilde{J}^N$, as well as the differences between the three costs, decrease with $N$, illustrating the convergence results of Section \ref{section: good aproximation}. Lastly, we observe that along the agent-to-destination proportion axis, the sign of $J^N(P)-\tilde{J}^N(P)$ varies. Indeed, when applying the continuum strategy for a finite population of size $N$, the actual agent-to-destination proportion is equal to $F^N(P)$ (See \eqref{eq: F^N}). Consequently, $\tilde{J}^N(P)$ is the cost incurred for the agent-to-destination proportion $F^N(P)$, which may slightly differ from $P$. Therefore, in Fig.\ref{fig: social costs}, we observe that the sign of $J^N(P)-\tilde{J}^N(P)$ is correlated with the relative variation of the social cost $J(P)$ at the agent-to-destination proportions $P$ and $F^N(P)$. However, as the population size increases and $F^N(P)$ approaches $P$, such discrepancy diminishes.

\section{Conclusion}
\label{section: conclusion}
This paper studies a linear quadratic collective discrete choice model.
We account for congestion by introducing a negative term in the running cost.
By exploiting the model's symmetries, we derive the agents'
optimal control strategy. 
Since computing such strategy may be computationally expensive, 
we develop an $\varepsilon-$optimal decentralized control strategy, where $\varepsilon$ tends 
to zero as the population size increases. 
This strategy significantly reduces the required computation compared to its optimal counterpart, 
as well as the amount of inter-agent information exchange.
In practice, the model may be useful for task assignment or resource allocation in multi-agent systems.
It is also of interest for future work to extend our analysis to agents with stochastic dynamics 
and to examine the case in which each agent can choose its exit time within the time interval $[0,T]$.
 
\bibliographystyle{ieeetr}
\bibliography{references}
\appendix

\section{Appendix}
\subsection{ Proof of Lemma \ref{lemma: social cost augmented state} }
\label{proof: social cost manipulation}

To prove the lemma, we note that
\begin{align*}
&\sum_{i=1}^N [x_i- \bar{x}^N]_{R_x}^2 \\
&= [X]^2_{R_x^{(N)}} - 2 (\bar{x}^N)^T R_x \left( \sum_{i=1}^N x_i \right) + N [\bar{x}^N]_{R_x}^2 \\
&= [X]^2_{R_x^{(N)}} - N [\bar{x}^N]_{R_x}^2.
\end{align*}
Moreover,
\begin{align*}
N [\bar{x}^N]_{R_x}^2 &= \frac{1}{N} \left( \sum_{i=1}^N x_i \right)^T R_x \left( \sum_{i=1}^N x_i \right) \\
&= \frac{1}{N} ((1_N^T \otimes I_n) X)^T R_x ((1_N^T \otimes I_n) X) \\
&= \frac{1}{N} X^T (1_N \otimes I_n) R_x (1_N \otimes I_n)^T X 
\end{align*}
\begin{align*}
\text{so }\quad  N [\bar{x}^N]_{R_x}^2 &= \frac{1}{N} X^T \left( 1_{N N} \otimes R_x \right) X.
\end{align*}
Hence,
\begin{align*}
&-\sum_{i=1}^N [x_i- \bar{x}^N]_{R_x}^2 + \sum_{i=1}^N [x_i - d_{\lambda_i}]^2_{R_d} 
\\
&= [X]^2_{-R_x^{(N)} + \frac{1}{N} \left( 1_{N N} \otimes R_x \right) + R_d^{(N)} } - 2 d_{\lambda^N}^T R_d^{(N)} X + [{d_{\lambda^N}}]^2_{R_d^{(N)}},
\end{align*}
and the rest of the argument leading to \eqref{eq:social cost augmented state} is straightforward.

\subsection{Proof of the Lemma \ref{lemma: dynamic programming} }
\label{proof: dynamic programming}
Let $\mathcal{H}^N$ and $V^N$ denote the Hamiltonian and Value Function associated with the cost $J^N_{\text{soc}}$ in (4). Specifically, we have
\begin{align}
    \mathcal{H}^N &= \frac{1}{2N} \left(
    [X]^2_{Q^N} - 2 d_{\lambda^N}^\top R_d^{(N)} X 
    + [d_{\lambda^N}]^2_{R_d^{(N)}} + [U]^2_{R_u^{(N)}} \right) \nonumber \\
    &\quad + (\nabla_X V^N)^\top (A^{(N)} X + B^{(N)} U). \nonumber
\end{align}

From the Hamilton-Jacobi-Bellman equation, we have:
\begin{align}
\label{HJB}
    &\frac{\partial V^N(t,X)}{\partial t} + \min_U \left\{ \mathcal{H}^N \right\} = 0, \nonumber \\
    & \text{with } V^N(T,X) = \frac{1}{2N} [X(T) - d_{\lambda^N}]^2_{M^{(N)}}.
\end{align}

We seek a quadratic form for the value function, given by:
\begin{equation}
\label{eq:value-function}
    V^N(t,X) = \frac{1}{2N} X^\top \varphi^N(t) X 
    + \frac{1}{N} (\Psi^N(t))^\top X + \chi^N(t), \nonumber
\end{equation}
where $\varphi^N(t) \in \mathbb{R}^{Nn \times Nn}$, $\chi^N(t) \in \mathbb{R}$, and
\begin{align}
\Psi^N(t) = \operatorname{col}(\Psi^N_1(t), \ldots, \Psi^N_N(t)) \in \mathbb{R}^{Nn}, \nonumber \\
\Psi^N_i(t) \in \mathbb{R}^n \;\text{for all } i \in \{1, \ldots, N\}. \nonumber
\end{align}

Since $R_u^{(N)}$ is positive definite, $\mathcal{H}^N$ is strictly convex with respect to $U$. Therefore, the optimal control, if it exists, is unique. Setting 
\(\frac{\partial \mathcal{H}^N}{\partial U} = 0\), we find that the optimal control $U^*$ satisfies:
\begin{align}
    U^* = -\left(R_u^{(N)}\right)^{-1}(B^{(N)})^\top \left( \varphi^N X + \Psi^N \right). \nonumber
\end{align}

Substituting $U = U^*$ in \eqref{HJB}, we obtain:
\begin{align}
    \frac{1}{2} X^\top & \frac{d\varphi^N}{dt} X 
    + \left( \frac{d \Psi^N}{dt} \right)^\top X + N \frac{d\chi^N}{dt} \nonumber \\
    &= \frac{1}{2} X^\top \left(-Q^N - \varphi^N A^{(N)} 
    - (A^{(N)})^\top \varphi^N \right. \nonumber \\
    &\quad \left. + \varphi^N B^{(N)} (R_u^{(N)})^{-1} (B^{(N)})^\top \varphi^N \right) X \nonumber \\
    &\quad + \left( (d_{\lambda^N})^\top R_d^{(N)} 
    - (\Psi^N)^\top A^{(N)} \right. \nonumber \\
    &\quad \left. + (\Psi^N)^\top B^{(N)} (R_u^{(N)})^{-1} (B^{(N)})^\top \varphi^N \right) X \nonumber \\
    &\quad - \frac{1}{2} (d_{\lambda^N})^\top R_d^{(N)} d_{\lambda^N} \nonumber \\
    & \quad + \frac{1}{2} (\Psi^N)^\top B^{(N)} (R_u^{(N)})^{-1} (B^{(N)})^\top \Psi^N. \nonumber
\end{align}

Thus, $\varphi^N$, $\Psi^N$, and $\chi^N$ must satisfy:
\begin{align}
    \frac{d \varphi^N}{dt} &= (\varphi^N)^\top S^{(N)} \varphi^N 
    - \varphi^N A^{(N)} - (A^{(N)})^\top \varphi^N - Q^N, \nonumber
 \\
    \frac{d \Psi^N}{dt} &= \left( \varphi^N S^{(N)} - (A^{(N)})^\top \right) \Psi^N 
    + R_d^{(N)} d_{\lambda^N}, \nonumber
 \\
    \frac{d \chi^N}{dt} &= \frac{1}{2N} \left( [\Psi^N]^2_{S^{(N)}} 
    - [d_{\lambda^N}]^2_{R_d^{(N)}} \right), \nonumber
\end{align}
with terminal conditions:
\begin{align}
    \varphi^N(T) &= M^{(N)}, \quad 
\Psi^N(T) = -M^{(N)} d_{\lambda^N}, \nonumber\\
\chi^N(T) &= \frac{1}{2N} [d_{\lambda^N}]^2_{M^{(N)}}. \nonumber
\end{align}
\subsection{ Proof of Lemma \ref{lemma: varphi decomposition} }
\label{proof: varphi decomposition}

Denote by $\pi(i,j) \in \mathbb{R}^{N \times N}$ the permutation matrix exchanging the $i^{th}$ 
and $j^{th}$ row of the identity matrix $I_N$, and let $\Pi(i,j)=\pi(i,j) \otimes I_n$. 
For $1 \leq i,j \leq N$ denote by $\varphi^N_{ij}$ the $(i,j)$-block of size 
$n \times n$ of the matrix $\varphi^N \in \mathbb R^{Nn\times Nn}$, solution of \eqref{eq:riccati all population}. 
Note that for any matrix $M \in \mathbb R^{N\times N}$, $\pi(i,j) M$ exchanges 
the rows $i$ and $j$ of $M$ and $M \pi(i,j)$ exchanges the columns $i$ and $j$ of $M$.
From this, one can see that the matrices $S^{(N)}$, $A^{(N)}$, $Q^N$ and $M^{(N)}$ appearing 
in \eqref{eq:riccati all population} and \eqref{eq:phi Psi terminal} are left unchanged 
when multiplied by $\Pi(i,j)$ on the left and on the right. Moreover, note that 
$\Pi(i,j)=\Pi(i,j)^T=\Pi(i,j)^{-1}$. Hence, if $\varphi^N$ is a solution of \eqref{eq:riccati all population},
so is $\Pi(i,j) \varphi^N \Pi(i,j)$, for any $1 \leq i, j \leq N$, and by unicity of the solution 
for the same terminal condition, all these matrices must be equal. This immediately implies that 
the diagonal blocks $\varphi^N_{ii}$ are all
equal, for $1 \leq i \leq N$, and the off-diagonal blocks $\varphi^N_{ij}$ for $j \neq i$, $1 \leq i,j \leq N$,
are also all equal. In other words, we can decompose $\varphi^N$ as in \eqref{eq:varphi decomposition}.
Finally, rewriting \eqref{eq:riccati all population} by block, we get
\eqref{eq:diagonal component for the riccati} and \eqref{eq:off-diagonal component for the riccati}.

\subsection{ Proof of Lemma \ref{lemma: psi decomposition}}
\label{proof: psi decomposition}

From \eqref{eq:big psi}, \eqref{eq:psi all population} and \eqref{eq:phi Psi terminal}, we have that 
for any $i \in \iset{N}$,
\begin{align}
    \frac{d}{d t} \Psi^N_i &= \left( \varphi^N_1S-A^T-\frac{1}{N}\varphi^N_2 S \right) \Psi^N_i+R_d d_{\lambda_i} \nonumber \\
    &+\frac{1}{N}\sum\limits_{j=1}^{N}\varphi^N_2 S \Psi^N_j, \;\; \Psi^N_i(T)=-Md_{\lambda_i}. \label{eq: sum by agent}
\end{align} 
Hence, for any two agents $i, j \in \iset{N}$, we have
\begin{align*}
    &\frac{d}{d t} (\Psi^N_i-\Psi^N_j)= \left( \varphi^N_1S-A^T-\frac{1}{N}\varphi^N_2S \right) (\Psi^N_i-\Psi^N_j) \\
    & \;\; +R_d (d_{\lambda_i}- d_{\lambda_j}), \;\; (\Psi^N_i-\Psi^N_j)(T)=-M(d_{\lambda_i}-d_{\lambda_j}). 
\end{align*}
Therefore, if $\lambda_i=\lambda_j$, i.e., the two agents have the same destination, 
then $0$ is a solution to this last differential equation and by unicity of the
solution we get $\Psi^N_i=\Psi^N_j$. 
We then obtain the conclusion of Lemma \ref{lemma: psi decomposition}, 
and \eqref{eq:psi principle components} is obtained by grouping the terms of the sum in \eqref{eq: sum by agent}
by destination.

\subsection{ Proof of Lemma \ref{lemma: psi form finite}}
\label{proof: psi form finite}

Let's assume that the solution of \eqref{eq:psi principle components} is of the form \eqref{eq: linear psi}.
Let $\bar{\varphi}_1^N=\varphi^N_1S-A^T-\frac{1}{N}\varphi^N_2S$, and $\bar{\varphi}_2^N=\varphi^N_2S$.
Injecting \eqref{eq: linear psi} in \eqref{eq:psi principle components} and recalling 
$P_D^N=1-\sum_{j=1}^{D-1} P_j^N$, we obtain
\begin{align*}
\frac{d}{d t}\psi^N_j &= \sum_{k=1}^{D-1} P^N_k \frac{d \alpha^N_k}{dt} - \frac{d \beta^N_j}{dt} \\
&= \bar{\varphi}_1^N \left[ \sum_{k=1}^{D-1}\alpha^N_k P^N_k -\beta^N_j \right] \\
& \quad + \sum\limits_{k=1}^{D} P_k^N \bar{\varphi}_2^N \left[ \sum_{l=1}^{D-1} \alpha^N_l P^N_l -\beta^N_k \right] + R_d d_j \\
\end{align*}
\begin{align*}
\frac{d}{d t}\psi^N_j &=\bar{\varphi}_1^N \left[ \sum_{k=1}^{D-1} \alpha^N_k P^N_k -\beta^N_j \right] + R_d d_j \\
& \quad +\sum\limits_{k=1}^{D-1}P_k^N \bar{\varphi}_2^N \left[ \sum_{l=1}^{D-1} \alpha^N_l P^N_l -\beta^N_k \right] \\
& \quad + \left( 1-\sum\limits_{k=1}^{D-1}P_k^N \right) \bar{\varphi}_2^N\left[ \sum_{l=1}^{D-1} \alpha^N_l P^N_l -\beta^N_D \right] \\
\end{align*}
so finally
\begin{align*}
\frac{d}{d t}\psi^N_j &= \sum_{k=1}^{D-1} P^N_k  \left[ (\bar{\varphi}_1^N+\bar{\varphi}_2^N) \alpha^N_k 
+ \bar{\varphi}_2^N(\beta_D^N - \beta_k^N) \right]  \\
& \quad -\bar{\varphi}_1^N \beta^N_j -\bar{\varphi}_2^N\beta^N_D + R_d d_j, 
\end{align*} 
with $\sum_{k=1}^{D-1} \alpha^N_k(T) P^N_k +\beta^N_j(T)=M d_j$ as terminal condition.
By identifying the terms, we get that the coefficients $\alpha^N_j$ and $\beta^N_j$ should 
respectively satisfy the linear differential equations \eqref{eq: finite alpha differential equation} 
and \eqref{eq: finite beta differential equation}, for which a unique solution exists.
Consequently, the solution candidate \eqref{eq: linear psi} satisfies the relation \eqref{eq:psi principle components}. 
By unicity of the solution to the linear differential equation \eqref{eq:psi principle components}, 
we conclude that \eqref{eq: linear psi} is the solution for \eqref{eq:psi principle components}.

\subsection{ Proof of Lemma \ref{lemma: chi form finite}}
\label{proof: chi form finite}

Note that, from \eqref{eq: linear psi} and the definition \eqref{eq: alphaN matrix}, we have
$\psi_j^N(t) = \alpha^N(t) P^N - \beta^N_j(t)$.
We then use the result of Lemma \ref{lemma: psi decomposition} and inject \eqref{eq: linear psi} 
in \eqref{eq:chi bar} to get

\begin{align*}
&\frac{d}{d t}\chi^N =\frac{1}{2} \sum\limits_{j=1}^{D} P_j^N \left( [\psi_j^N]^2_{S}-[d_j]^2_{R_d} \right)
\end{align*}
\begin{align*}
&\frac{d}{d t}\chi^N 
=\frac{1}{2} \sum\limits_{j=1}^{D} P_{j}^N \left( \left[ \alpha^N P^N -\beta_j^N \right]^2_{S}-[d_j]^2_{R_d}\right) \\
    &= \frac{1}{2} (P^N)^T(\alpha^N)^T S \alpha^NP^N-\sum_{j=1}^{D} P^N_j(\beta^N_j)^T S \alpha^NP^N   \\
    &+\frac{1}{2}\sum_{j=1}^{D}P_j^N\left([\beta_j^N]^2_{S}-[d_j]^2_{R_d}\right)
\end{align*}
\begin{align*}
&\frac{d}{d t}\chi^N = (P^N)^T W^N P^N +\frac{1}{2}\sum_{j=1}^{D}P_j^N\left([\beta_j^N]^2_{S}-[d_j]^2_{R_d}\right),
\end{align*}
with $W^N$ defined in \eqref{eq: W N}.

\subsection{Proof of Lemma \ref{lemma:value function rewritten}}
\label{proof:value function rewritten}

Using \eqref{eq: value function} and \eqref{eq: linear psi}, we have that
\begin{align*}
&V^N(t,\lambda^N)
= \frac{1}{2N}X^T\varphi^N X+\frac{1}{N}(\Psi^N)^T X+\chi^N \\
&= \frac{1}{2N}\left(\sum_{i=1}^{N}x_i^T \varphi_1^N x_i
+\sum_{i=1}^{N} \sum_{\substack{j=1\\j \neq i}}^{N} \frac{1}{N} x_i^T \varphi_2^N x_j \right) \\
&+\frac{1}{N} \sum\limits_{i=1}^{N} \left( \left( \sum_{j=1}^{D-1} P^N_j \alpha^N_j \right) 
-\beta^N_{\lambda_i}\right)^T x_i+\chi^N(t) 
\end{align*}
\begin{align*}
&=\frac{1}{2N}\sum_{i=1}^{N}x_i^T\varphi_1^N x_i
+\frac{1}{2} \left( \frac{1}{N}\sum_{i=1}^{N}x_i \right)^T \varphi_2^N \left( \frac{1}{N}\sum_{i=1}^{N} x_i \right)\\ 
&-\frac{1}{2N^2}\sum_{i=1}^{N}x_i^T\varphi_2^N x_i+\frac{1}{N}\sum\limits_{i=1}^{N}\left(\sum\limits_{j=1}^{D-1}P^N_j(\alpha_j^N)^Tx_{i} \right.\\
&\left.+\frac{1}{2} \left( |x_{i}-\beta_{\lambda_i}^N|_2^2-|x_i|_2^2- |\beta_{\lambda_i}^N|_2^2 \right) \right)+\chi^N(t)
\end{align*}
\begin{align*}
V^N(t,\lambda^N) &=\frac{1}{2N}\sum_{i=1}^{N} x_i^T \left( \varphi_1^N-\frac{\varphi_2^N}{N} \right) x_i
+\frac{1}{2}(\bar{x}^N)^T\varphi_2^N \bar{x}^N\\ 
&+\sum\limits_{j=1}^{D-1}P^N_j(\alpha_j^N)^T\bar{x}^N-\frac{1}{2}\sum\limits_{j=1}^{D}P^N_j|\beta_{j}^N|_2^2\\
&+\frac{1}{2N}\sum\limits_{i}=1^{N}|x_{i}-\beta_{\lambda_i}^N|_2^2-\frac{1}{2N}\sum\limits_{i=1}^{N}|x_i|_2^2+\chi^N(t).
\end{align*}

\subsection{Proof of Theorem \ref{thm: bound on time horizon}}
\label{proof:bounding escape time}

Define the matrix $Q_{\min}^N=I_N \otimes (R_d-R_x)$, so that

\begin{align*}
Q^N-Q_{\min}^N&=\frac{1}{N} 1_{NN} \otimes R_x. 
\end{align*} 
The eigenvalues of the Kronecker product of any two square matrices are the
pairwise products of the eigenvalues of both matrices. 
The only non-zero eigenvalue of the rank $1$ matrix $\frac{1}{N} 1_{NN}$ is $1$. 
Since $R_x$ is positive definite, we conclude that 
$Q_{\min}^N\preccurlyeq Q^N$.
Therefore, by using the order preserving property for 
Hermitian Riccati equations \cite[Theorem 4.1.4]{abou2012matrix}, we conclude that
$\Delta_{esc}(Q_{\min}^N) \leq \Delta_{esc}(Q^N)$, where, to simplify the notation, 
we have redefined $\Delta_{esc}(\tilde Q) = \Delta_{esc}(S^{(N)},A^{(N)},\tilde Q,M^{(N)})$,
for $\tilde Q \in \mathbb R^{Nn\times Nn}$ symmetric.

Suppose now that $T < \Delta_{esc}(Q_{\min}^N)$. Then
\[
T - \Delta_{esc}(Q^N) \leq T - \Delta_{esc}(Q_{\min}^N) < 0,
\] 
so the Riccati equation \eqref{eq:riccati all population} has a solution 
over the time interval $[0,T]$.
Finally, since all the matrices 
in $\textit{HDRE}(T,S^{(N)},A^{(N)},Q_{\min}^N,M^{(N)})$
are block diagonal  
with repeated diagonal blocks, we have that
$\Delta_{esc}(Q_{\min}^N)=\Delta_{esc}(S,A,R_d-R_x,M)$. This gives the sufficient condition
of the theorem.

\subsection{Proof of Lemma \ref{lemma: dual optimum}}
\label{proof: dual function lemma}

Define $c(x,j) = |x - \beta_j(0)|^2$, for any $x$ in $E$ and $j \in \iset{D}$.
By Kantorovitch duality \cite[Theorem 5.10]{Villani2009OptimalNew}, we can write for the OT 
problem \eqref{eq:semi-discrete kantorovich prime} 
\begin{equation}    \label{eq: dual MKP}
C(P) = \max_{(\phi,g) \in \Phi^c} \int_E \phi(x) d \mathcal P_0(x) + \sum_{j=1}^D P_j g_j,
\end{equation}
where, denoting $\overline{\mathbb R} = \mathbb R \cup \{-\infty\}$, we have
\begin{align}  
& \Phi^c = \Big \{ \phi: \mathbb R^n \to \overline{\mathbb R}, g \in \overline{\mathbb R}^D : \int_E |\phi(x)| d \mathcal P_0(x) < \infty, \label{eq: condition on dual potentials} \\
& \sum_{j=1}^D |g_j| P_j < \infty; \phi(x)+g_j \leq c(x,j), \forall x \in E, \forall j \in \iset{D}  \Big \}.    \nonumber 
\end{align}
The notation $\max$ in \eqref{eq: dual MKP} means that an maximizing pair $(\tilde \phi,\tilde g)$ in $\Phi^c$ exists.
However, the summability condition on $g$ in \eqref{eq: condition on dual potentials} does not rule out 
a priori that $\tilde g_j = - \infty$ when $P_j = 0$, see Remark \ref{remark: finite dual optimum}.
Nonetheless, due to the inequality constraint in $\Phi^c$, starting from $\tilde \phi$, $\tilde g$, 
we can replace $\tilde \phi$ by its $c$-transform \cite{Villani2009OptimalNew}, defined for all $x$ in $E$ as
\begin{equation}    \label{eq: first c transform}
\phi^*(x) \coloneqq \min_{j \in \iset{D}} \{ c(x,j) - g^*_j \} = \min_{j: g^*_j > -\infty} \{ c(x,j) - g^*_j \},
\end{equation}
without reducing the cost in \eqref{eq: dual MKP}. 
In \eqref{eq: first c transform} the elements $g_j^* = -\infty$ lead to $c(x,j) - g_j^*$ and hence 
can be removed from the minimization. Note that this leads to the form \eqref{eq:dual} for the dual function,
but more importantly here, $\phi^*$ is continuous because $c(\cdot,j)$ is continuous for each $j \in \iset{D}$.
Similarly, we can then replace $\tilde g$ by its $c$-transform $g^*$, with components
\begin{equation}    \label{eq: second c-transform}
g^*_j = \inf_{x \in E} \{ c(x,j) - \phi^*(x) \}, \;\; j \in \iset{D},
\end{equation}
without reducing the cost in \eqref{eq: dual MKP}. But since $\phi^*$ is continuous, $c(\cdot,j)$
is continuous, and $E$ is compact by Assumption \ref{asspt: E compact}, the infimum is attained
in \eqref{eq: second c-transform} and $g^*_j$ is finite, for each $j \in \iset{D}$. This proves
the second part of the lemma, i.e., $g^* \in \mathbb R^D$ is an maximizer of the dual function.

For the first part, $g \mapsto J^D(g,P)$ is concave as a general property of dual functions, or more 
directly one can observe that $\min_j\{c(x,j) - g_j\}$ is concave in $g$ as the minimum of linear functions,
and concavity is preserved by integration.

\subsection{Proof of Lemma \ref{prop: psi cost part is convex}}
\label{proof: convexity of C(P)}

Recall first the definition of the Wasserstein $2$-distance between two measures $\mu$ and $\nu$ on $\mathbb R^n$
\[
W_2(\mu,\nu) = \left( \inf_{\gamma \in \Gamma(\mu,\nu)} \int_{\mathbb R^n \times \mathbb R^n} |x-y|_2^2 \; d\gamma(x,y) \right)^{1/2},
\]
where $\Gamma(\mu,\nu)$ is the set of all joint measures in $\mathcal M(\mathbb R^n \times \mathbb R^n)$ with
marginals $\mu$ and $\nu$, i.e., 
\[
\gamma(A \times \mathbb R^n) = \mu(A), \;\; \gamma(\mathbb R^n \times B) = \nu(B),
\]
for all measurable sets $A$ and $B$.

Define, for any $P \in \Omega$, the probability measure $\mathcal P_T(P) = \sum_{i=1}^D P_i \delta_{\beta_i(0)}$
on $\mathbb R^n$, where $\delta_x$ denotes the Dirac measure at $x$.
Notice from \eqref{eq:semi-discrete kantorovich prime} that 
\begin{align*}
C\left(P\right)=W_2\left( \mathcal{P}_0,\mathcal{P}_T(P) \right)^2.
\end{align*}
Consider a sequence $\{P^{(n)}\}_{n \geq 0}$ in $\Omega$ converging to some vector $P \in \Omega$. 
Then $\mathcal P_T(P^{(n)})$ converges weakly to $\mathcal P_T(P)$, so by \cite[Theorem 6.9]{Villani2009OptimalNew},
under Assumption \ref{asspt: E compact}, 
$\lim_{n \to \infty} C(P^{(n)}) \to C(P)$, which shows the continuity of $C$ on $\Omega$.

For the convexity of $C$, consider two vectors $P^1,P^2$ in $\Omega$ and $\lambda \in [0,1]$. 
The cost function $c(x,y) = |x-y|_2^2$ is continuous and non-negative on $\mathbb R^n \times \mathbb R^n$.
Define the discrete measure $\Lambda$ on the set $\{1,2\}$ by
$\Lambda(\{1\}) = \lambda$, $\Lambda(\{2\}) = 1-\lambda$.

According to \cite[Theorem 4.8]{Villani2009OptimalNew}, we then have
\[
C(\lambda P^1+(1-\lambda) P^2) \leq \lambda C( P^1)+(1-\lambda) C(P^2),
\]
i.e., $C$ is convex on $\Omega$.

 \subsection{Proof of Lemma \ref{prop: chi is convex}}
\label{proof: convexity of chi}

By integrating \eqref{eq: chi polynomial infinite}, we have
\begin{align*}
\chi(0,P) &= \chi(T,P) - P^T \left(\int_0^T W(\tau) d\tau \right) P \\
&- \frac{1}{2} \sum_{j=1}^{D} P_j \int_{0}^T  [\beta_j(\tau)]^2_{S}-[d_j]^2_{R_d} d \tau.
\end{align*}
Since $\chi(T,P)$ is linear in $P$, to prove the convexity of $\chi(0,P)$ 
on $\Omega$ it is sufficient to prove that $P \to P^T W(t) P$
is concave on $\Omega$, for all $t$. Let us parametrize $\Omega$ as
\begin{align*}
\tilde \Omega = \{ &\begin{bmatrix} P_1  \ldots P_{D-1} \end{bmatrix}^T \in \mathbb R^{D-1} \Big | \\
& P_i \geq 0, \forall i \in \iset{D-1}, \sum_{i=1}^{D-1} P_i \leq 1 \},
\end{align*}
so $P \in \Omega$ if and only if $P = \begin{bmatrix} \tilde P^T, (1-1_{D-1}^T \tilde P) \end{bmatrix}^T$
for $\tilde P \in \tilde \Omega$. 
Let $\tilde \alpha = \begin{bmatrix} \alpha_1 \ldots \alpha_{D-1} \end{bmatrix}$,
$\tilde \beta = \begin{bmatrix} \beta_1 \ldots \beta_{D-1} \end{bmatrix}$.
Note that 
\[
W = \begin{bmatrix} \tilde W & 0_{(D-1) \times 1} \\ -\beta_D^T S \tilde \alpha & 0\end{bmatrix},
\;\; \tilde W = \left( \frac{\tilde \alpha}{2}- \tilde \beta \right)^T S \tilde \alpha,
\]
and, for $P = \begin{bmatrix} \tilde P^T, (1-1_{D-1}^T \tilde P) \end{bmatrix}^T$,
\[
P^T W P = \tilde P \left( \tilde W + 1_{D-1} \beta_D^T S \tilde \alpha \right) \tilde P 
- \beta_D^T S \tilde \alpha \tilde P.
\]
As a result, the function is $P \to P^T W P$ is concave on $\Omega$ if the symmetric part
of the matrix 
\begin{equation}    \label{eq: W tilde +}
\tilde W + 1_{D-1} \beta_D^T S \tilde \alpha = 
\left(\tilde \alpha - \left( \tilde \beta - \beta_D 1_{D-1}^T \right) \right)^T S \tilde \alpha
- \frac{1}{2} \tilde \alpha^T S \tilde \alpha
\end{equation}

is negative semi-definite.

Let $d=[d_1-d_D,\hdots,d_{D-1}-d_D] \in \mathbb{R}^{n \times (D-1)}$. 
Let $\bar{\varphi}_1=\varphi_1S-A^T$, $\bar{\varphi}_2=\varphi_2S$ and
let $\Phi_{\bar{\varphi}_1}$ and $\Phi_{\bar{\varphi}_1+\bar{\varphi}_2}$ be
the state transition matrices of $\bar{\varphi}_1$ and $\bar{\varphi}_1+\bar{\varphi}_2$ respectively.
From \eqref{eq: beta differential equation}, we have for all $j \in \iset{D-1}$
\begin{align}
\frac{d}{dt} (\beta_j-\beta_D) &= \bar{\varphi}_1 (\beta_j-\beta_D) - R_d (d_j-d_D),  \nonumber \\
(\beta_j-\beta_D)(T) &= M(d_j-d_D). \label{eq: beta new}
\end{align} 
Therefore,

\begin{align}
\tilde{\beta}(t) - \beta_D(t) 1_{D-1}^T = \bar{\phi}_{\beta}(t) \, d,  \label{eq: new beta} 
\end{align} 
with 
\[
\bar{\phi}_{\beta}(t) = \int_t^T \Phi_{\bar{\varphi}_1}(t,\tau)R_d \, d \tau+\Phi_{\bar{\varphi}_1}(t,T) M.
\]
By subtracting \eqref{eq: beta new} from \eqref{eq: alpha differential equation}, 
we get that for all $j \in \iset{D-1}$,
\begin{align*}
&\frac{d}{dt} \left(\alpha_j-(\beta_j-\beta_D)\right)=(\bar{\varphi}_1+\bar{\varphi}_2) 
\left(\alpha_j-(\beta_j-\beta_D)\right)\\
& \;\; + R_d (d_j-d_D) \;, \; \left(\alpha_j-(\beta_j-\beta_D)\right)(T)=-M(d_j-d_D). \nonumber
\end{align*}
Hence,
\begin{align*}
\tilde \alpha(t) - \left( \tilde \beta(t) - \beta_D(t) 1_{D-1}^T \right)  
&= -\bar{\phi}_{\alpha-\beta}(t) \, d,
\end{align*}
with
\[
\bar{\phi}_{\alpha-\beta}(t) = \int_t^T \Phi_{\bar{\varphi}_1+\bar{\varphi}_2}(t,\tau)R_d d\tau + \Phi_{\bar{\varphi}_1+\bar{\varphi}_2}(t,T) M.
\]
From \eqref{eq: alpha differential equation} and \eqref{eq: new beta},
we establish the following
\begin{align*}
\tilde{\alpha}(t) &=
-\int_T^t \Phi_{\bar{\varphi}_1+\bar{\varphi}_2}(t,\tau) \, \bar{\varphi}_2(\tau) 
\left(\tilde{\beta} - \beta_D 1_{D-1}^T\right)(\tau) d \tau  \\
&= \left( \int_t^T \Phi_{\bar{\varphi}_1+\bar{\varphi}_2}(t,\tau) \bar{\varphi}_2(\tau) 
\bar{\phi_{\beta}}(\tau) d\tau \right) d
=:\bar{\phi}_{\alpha}(t) \, d.
\end{align*}

Let us call a diagonal matrix with all nonnegative or all positive diagonal coefficients 
a diagonal nonnegative or positive matrix respectively.
Under Assumption \ref{assmp: diagonal matrices}, $S$ is a diagonal 
nonnegative matrix. Furthermore, under Assumption \ref{assmp: diagonal matrices} and using 
\cite[Theorem 4.1.6]{abou2012matrix}, we can show that the solution $\varphi_2$ of 
the Riccati equation \eqref{eq:off-diagonal component for the riccati inf} is a diagonal nonnegative 
matrix. 
Therefore, 
$\bar{\varphi}_2$
is a diagonal nonnegative matrix.
Besides, any state transition matrix of a diagonal matrix is a diagonal positive matrix. Therefore, since the matrices $A,B,R_X,R_d,R_u,M$  
are diagonal, 
the matrices $\bar{\phi}_{\alpha-\beta}$, $\bar{\phi}_{\beta}$ and $\bar{\phi}_{\alpha}$ are diagonal nonnegative matrices.
As a result, $\bar{\phi}_{\alpha-\beta}^T S \bar{\phi}_{\alpha}$ is a diagonal
nonnegative matrix and so, from \eqref{eq: W tilde +},
\begin{align*}\tilde W(t) + 1_{D-1} \beta_D^T(t) S \tilde \alpha(t)
&=-d^T \bar{\phi}_{\alpha-\beta}^T S \bar{\phi}_{\alpha} d-\frac{1}{2}\tilde{\alpha}^TS\tilde{\alpha}
\end{align*} 
is a negative semi-definite matrix, for all $t \in [0,T]$. 
Hence, the function $P \rightarrow \chi(0,P)$ is convex under Assumption \ref{assmp: diagonal matrices}.
\subsection{Proof of Lemma \ref{lemma: value function uniform convergence}}
\label{proof: value function uniform convergence}

For any $P \in \Omega$, we have from \eqref{eq: J bar N} and \eqref{eq:simplified value function infinite}
\begin{align*}
&\bar J^N(P) - J(P) = \\
&\frac{1}{2} \left( \Delta^N_I(P) + \Delta^N_{II}(P) + \Delta^N_{III} + \Delta^N_{IV} + \Delta^N_{V} + \Delta_{VI}\right),
\end{align*}
with
\[
\Delta^N_I(P) = \bar{C}^N(P)-C(P), \;\; \Delta^N_{II}(P) = \chi^N(0,P)-\chi(0,P),
\]
\[
\Delta^N_{III} = \int_{x \in E}|x|^2d \mathcal{P}_0(x) - \frac{1}{N}\sum\limits_{i=1}^{N}|x_i(0)|^2,
\]
\begin{align*}
\Delta^N_{IV} &= \frac{1}{N} \sum_{i=1}^{N}x_i(0)^T \left( \varphi^N_1(0)-\frac{\varphi^N_2(0)}{N} \right) x_i(0) \\
&- \int_{x \in E}x^T\varphi_1(0)x \, d \mathcal{P}_0(x),
\end{align*}
\[
\Delta^N_{V} = \bar{x}^N(0)^T \varphi^N_2(0) \bar{x}^N(0) - \bar{x}(0)^T \varphi_2(0)\bar{x}(0),
\]
\begin{align*}
\Delta^N_{VI} &= \sum_{j=1}^{D-1} 2 P_j \left( \alpha^N_j(0)^T \bar{x}^N(0)-\alpha_j(0)^T\bar{x}(0) \right) \\
& + \sum_{j=1}^{D} P_j \left( |\beta_{j}(0)|_2^2-|\beta_{j}^N(0)|_2^2 \right).
\end{align*}

By the strong law of large numbers (SLLN), $\Delta^N_{III} \to 0$ as $N \to \infty$, almost surely. We also have the following results.

\begin{lemma}   \label{lemma: alpha, beta, phi convergence}
The functions $\varphi^N_1$, $\varphi^N_2$, $\alpha^N_j$, $j \in \llbracket D-1 \rrbracket$ and 
$\beta^N_j$, $j \in \llbracket D \rrbracket$, solutions to \eqref{eq:diagonal component for the riccati}, 
\eqref{eq:off-diagonal component for the riccati}, \eqref{eq:off-diagonal component for the riccati}, 
\eqref{eq: finite alpha differential equation}, converge uniformly on $[0,T]$ to
the solutions $\varphi_1$, $\varphi_2$, $\alpha_j$ and $\beta_j$ 
of \eqref{eq:diagonal component for the riccati inf}, \eqref{eq:off-diagonal component for the riccati inf},
\eqref{eq: alpha differential equation}, \eqref{eq: beta differential equation} respectively, as $N \to \infty$.
\end{lemma}
\begin{lemma}   \label{lemma: C(P) convergence}
Almost surely, the function $P \mapsto \Delta^N_I(P)$ converges uniformly to $0$ on $\Omega$ as $N \to \infty$.
\end{lemma}
\begin{lemma}
\label{lemma: squared uniform convergence}
The functions $W^N$ defined in \eqref{eq: W N} and $[\beta^N_j]_S^2$, 
$j \in \llbracket D \rrbracket$ converge uniformly on $[0,T]$ to $W$ defined in \eqref{eq: W inf} and $[\beta_j]_S^2$ respectively, as $N \to \infty$.
\end{lemma}
\begin{lemma}   \label{lemma: chi convergence}
The function $P \mapsto \Delta_{II}^N(P)$ converges uniformly to $0$ on $\Omega$ as $N \to \infty$.
\end{lemma}

Using Lemma \ref{lemma: alpha, beta, phi convergence} and the SLLN, it is straightforward to show that
almost surely, $\Delta^N_{IV}$, $\Delta^N_V$ and $\Delta_{VI}^N$ converge to $0$ as $N \to \infty$.
Hence, with Lemma \ref{lemma: C(P) convergence} and Lemma \ref{lemma: chi convergence}, we conclude 
that with probability one, $\bar J^N$ converges uniformly to $J$ on $\Omega$. 
This proves Lemma \ref{lemma: value function uniform convergence}.
The rest of the section contains the proofs of Lemmas \ref{lemma: alpha, beta, phi convergence}, \ref{lemma: C(P) convergence},
\ref{lemma: squared uniform convergence} and \ref{lemma: chi convergence}.

\subsubsection{Proof of Lemma \ref{lemma: alpha, beta, phi convergence}}

For each $N$, the functions $\varphi^N_1$, $\varphi^N_2$, $\alpha^N_j$, $j \in \llbracket D-1 \rrbracket$ 
and $\beta^N_j$, $j \in \llbracket D \rrbracket$, form the solution $t \to \xi^N(t) \in \mathbb{R}^{2n^2+(2D-1)n}$ 
to the system of ordinary differential equations (ODEs) \eqref{eq:diagonal component for the riccati}, 
\eqref{eq:off-diagonal component for the riccati}, 
\eqref{eq: finite alpha differential equation} and \eqref{eq: finite beta differential equation},
which  we can rewrite as $\dot \xi^N = f^N(\xi^N)$, with $\xi^N(T) = \xi_T$, for some function $f^N$.
Similarly, the functions $\varphi_1$, $\varphi_2$, $\alpha_j$, $j \in \llbracket D-1 \rrbracket$ and 
$\beta_j$, $j \in \llbracket D \rrbracket$ form the solution $t \to \xi(t)$ to the system of ODEs
\eqref{eq:diagonal component for the riccati inf}, \eqref{eq:off-diagonal component for the riccati inf},
\eqref{eq: alpha differential equation}, \eqref{eq: beta differential equation},
which we can rewrite as $\dot \xi = f(\xi)$, with $\xi(T) = \xi_T$, for some function $f$.
Note that these two systems have the same terminal conditions.
It is straightforward to see that as $N \to \infty$, $f^N$ converges uniformly 
on compacts to $f$. Moreover, the solution to the limit system $\dot \xi = f(\xi)$, $\xi(T) = \xi_T$,
is unique by the standard properties of Riccati and linear ODEs. Hence, by 
\cite[Theorem 3]{Coppel1965StabilityEquations}, we conclude that $t \to \xi^N(t)$ converges
uniformly on $[0,T]$ to $t \to \xi(t)$ and $N \to \infty$.

\subsubsection{Proof of Lemma \ref{lemma: C(P) convergence}}

Recall first the definitions given at the beginning of Appendix \ref{proof: convexity of C(P)}.
For any $P \in \Omega$ and any initial conditions $x_1(0),x_2(0),\ldots,$ sampled independently 
according to $\mathcal P_0$, define the probability measures 
\begin{align}    
\mathcal{P}^N_0 &= \sum_{i=1}^N \frac{1}{N} \delta_{x_i(0)},  \label{eq: empirical init distrib} \\
\mathcal{P}^N_T(P) &=\sum_{i=1}^D P_{i} \delta_{\beta^N_i(0)}, \text{ and }
\mathcal{P}_T(P) =\sum_{i=1}^D P_{i} \delta_{\beta_i(0)} \nonumber
\end{align} 
on $\mathbb R^n$.
To simplify the notation we sometimes omit the argument $P$ from 
$\mathcal{P}^N_T$ and $\mathcal{P}_T$ in the following.
With the definitions \eqref{eq: discrete Kantorovich prime} and \eqref{eq:semi-discrete kantorovich prime},
we have
\begin{align*}
\bar{C}^N\left(P\right)=W_2\left(\mathcal{P}^N_0,\mathcal{P}^N_T\right)^2 \text{ and }
C\left(P\right)=W_2\left(\mathcal{P}_0,\mathcal{P}_T\right)^2, 
\end{align*}
and so
\begin{align}
|\bar{C}^N\left(P\right)-C\left(P\right)| & =  
|W_2\left(\mathcal{P}^N_0,\mathcal{P}^N_T\right)-W_2\left(\mathcal{P}_0,\mathcal{P}_T\right)| \times \nonumber \\
& \left( W_2\left(\mathcal{P}^N_0,\mathcal{P}^N_T\right)+W_2\left(\mathcal{P}_0,\mathcal{P}_T\right) \right). 
\label{eq: basic Wasserstein relation}
\end{align}
In the following, we show that the right-hand side of \eqref{eq: basic Wasserstein relation} converges 
uniformly on $\Omega$ to $0$, almost surely, by considering each term in the product.

Since $W_2$ defines a distance between probability measures with finite second moment \cite[Chapter 6]{Villani2009OptimalNew}, 
by the triangle inequality
\begin{align*}
W_2\left(\mathcal{P}^N_0,\mathcal{P}^N_T\right) - &W_2\left(\mathcal{P}_0,\mathcal{P}_T\right) 
\leq \\
&W_2\left(\mathcal{P}^N_0,\mathcal{P}_0\right)+W_2\left(\mathcal{P}_T,\mathcal{P}^N_T\right).
\end{align*}
Similarly,
\begin{align*}
W_2\left(\mathcal{P}_0,\mathcal{P}_T\right) - &W_2\left(\mathcal{P}^N_0,\mathcal{P}^N_T\right) \leq \\ 
& W_2\left(\mathcal{P}_0,\mathcal{P}^N_0\right)+W_2\left(\mathcal{P}^N_T,\mathcal{P}_T\right).
\end{align*}
Hence, by the symmetry property of the distance,
\begin{align}
|W_2\left(\mathcal{P}^N_0,\mathcal{P}^N_T\right) - &W_2\left(\mathcal{P}_0,\mathcal{P}_T\right)| \leq 
\label{eq: triangle ineq. Wasserstein}
\\
& W_2\left(\mathcal{P}^N_0,\mathcal{P}_0\right)+W_2\left(\mathcal{P}^N_T,\mathcal{P}_T\right). \nonumber
\end{align}
By considering the transportation plan sending all the mass $P_i$ at $\beta^N_i(0)$ to $\beta_i(0)$, we
obtain that
\begin{align*}
W_2\left(\mathcal{P}^N_T,\mathcal{P}_T\right) &\leq \sqrt{\sum_{i=1}^D P_i|\beta^N_i(0)-\beta_i(0)|_2^2} \\
& \leq \sqrt{\sum_{i=1}^D|\beta^N_i(0)-\beta_i(0)|_2^2},
\end{align*}
using the fact that $P_i \leq 1$, for all $i \in \iset{D}$. 
Hence, from Lemma \ref{lemma: alpha, beta, phi convergence},
\begin{align}   \label{eq: unif. cv partial}
\lim_{N \to \infty} \sup_{P \in \Omega} W_2 \left( \mathcal{P}^N_T,\mathcal{P}_T \right) = 0.
\end{align}
The random measure $\mathcal P_0^N$ converges weakly to $\mathcal P_0$, almost surely \cite{varadarajan1958convergence}.
Hence, by \cite[Theorem 6.9]{Villani2009OptimalNew} and under Assumption \ref{asspt: measure zero},
\begin{align}   \label{eq: a.s. cv P0}
\lim_{N \to \infty} W_2(\mathcal P_0^N,\mathcal P_0) = 0, \; \text{almost surely}.
\end{align}
Combining \eqref{eq: a.s. cv P0} with \eqref{eq: unif. cv partial} in \eqref{eq: triangle ineq. Wasserstein} 
shows that almost surely,
\begin{align}   \label{eq: convergence first term}
\lim_{N \to \infty} \sup_{P \in \Omega} 
|W_2\left(\mathcal{P}^N_0,\mathcal{P}^N_T\right) - W_2\left(\mathcal{P}_0,\mathcal{P}_T\right)| = 0.
\end{align}

For the other term in the product \eqref{eq: basic Wasserstein relation}, by the triangle inequality again
\begin{align*}
&W_2\left(\mathcal{P}^N_0,\mathcal{P}^N_T\right)+W_2\left(\mathcal{P}_0,\mathcal{P}_T\right) \\
&\leq W_2\left(\mathcal{P}^N_0,\mathcal{P}_0\right) + 2 W_2\left(\mathcal{P}_0,\mathcal{P}_T\right)
+ W_2\left(\mathcal{P}_T,\mathcal{P}^N_T\right).
\nonumber \end{align*} 
By \eqref{eq: unif. cv partial} and \eqref{eq: a.s. cv P0}, almost surely, there exists some constant $c_1$ and
some integer $N_0$ such that for all $N \geq N_0$,
\[
W_2\left(\mathcal{P}^N_0,\mathcal{P}_0\right) + \sup_{P \in \Omega} W_2\left(\mathcal{P}_T,\mathcal{P}^N_T\right) \leq c_1.
\]
Finally, by \cite[Corollary 6.10]{Villani2009OptimalNew}, $P \mapsto W_2\left(\mathcal{P}_0,\mathcal{P}_T(P) \right)$
is a continuous function of $P$ on $\Omega$, which is compact, hence it is uniformly bounded on $\Omega$.
We conclude from  that almost surely, the term 
$W_2\left(\mathcal{P}^N_0,\mathcal{P}^N_T\right)+W_2\left(\mathcal{P}_0,\mathcal{P}_T\right)$ 
in \eqref{eq: basic Wasserstein relation} is uniformly bounded on $\Omega$ for $N$ large enough,
and hence with \eqref{eq: convergence first term} obtain the result of the lemma, i.e.,
\[
\lim_{N \to \infty} \sup_{P \in \Omega} |\bar C^N(P) - C(P)| = 0, \text{ almost surely}.
\]
\subsubsection{Proof of Lemma \ref{lemma: squared uniform convergence}}
The coefficients $\alpha_j$ and $\beta_j, j \in \llbracket D-1 \rrbracket$ are solution of the differential equations \eqref{eq: alpha differential equation}, \eqref{eq: beta differential equation} respectively. Therefore, they are bounded over the time interval $[0,T]$. By the uniform convergence of the coefficients $\alpha_j^N$ and $\beta_j^N$ to $\alpha_j$ and $\beta_j, j \in \llbracket D-1 \rrbracket$ with respect to $t \in [0,T]$, we establish that $W^N$ defined in \eqref{eq: W N}, $\alpha_j^N$, and $\beta_j^N, j \in \llbracket D-1 \rrbracket$ are uniformly bounded. Since the product of bounded uniformly convergent functions is  uniformly convergent \cite[Exercise 2]{rudin1976principles}, 
we conclude from Lemma \ref{lemma: alpha, beta, phi convergence} that the functions $W^N$ and $[\beta^N_j]_S^2$ converge uniformly on $[0,T]$ to $W$ defined in \eqref{eq: W inf} and $[\beta_j]_S^2$ respectively, with $j \in \llbracket D \rrbracket$.
\subsubsection{Proof of Lemma \ref{lemma: chi convergence}}
We need to show that $\chi^N(0,\cdot)$ converges uniformly to $\chi(0,\cdot)$ on $\Omega$,
with $\chi^N$ and $\chi$ the solutions of \eqref{eq: chi polynomial finite} and \eqref{eq: chi polynomial infinite}.
Let $\Delta^N_\chi(t,P) = \chi(t,P) - \chi^N(t,P)$. First, $\Delta^N_\chi(T,P) = 0$.

Then, denoting by $\|\cdot\|_2$ the induced $2$-norm,
\begin{align*}
& \left| \frac{\partial}{\partial t} \Delta^N_\chi(t,P) \right| \\
&= \left| P^T (W-W^N) P +\frac{1}{2}\sum_{j=1}^{D}P_j\left([\beta_j]^2_{S}-[\beta_j^N]^2_{S}\right) \right| \\
& \leq \| W-W^N \|_2 \; | P |_2^2 
+\frac{1}{2}\sum_{j=1}^{D}P_j \left| [\beta_j]^2_{S}-[\beta^N_j]^2_{S} \right| \\
&\leq \| W(t)-W^N(t)\|_2 
+\frac{1}{2}\sum_{j=1}^{D}|[\beta_j(t)]^2_{S}-[\beta^N_j(t)]^2_{S}|. 
\end{align*}
Therefore,
\begin{align}
|\Delta_\chi(0,P)| &\leq \int_0^T \| W(\tau)-W^N(\tau)\|_2d \tau  \nonumber\\
&+\frac{1}{2}\sum_{j=1}^{D}\int_0^T \left| [\beta_j(\tau)]^2_{S}-[\beta^N_j(\tau)]^2_{S} \right| d \tau. 
\label{eq: chi diff ineq}
\end{align}
The upper bound in \eqref{eq: chi diff ineq} does not depend on the parameter $P$, 
and using the result of Lemma \ref{lemma: squared uniform convergence}
we conclude  that $\chi^N(0,P)$ converges uniformly to $\chi(0,P)$ with respect to $P$. 
\subsection{Proof of Lemma \ref{lemma: conv optimal cost to limit}}
\label{proof: conv optimal cost to limit}

Define $P^N_{opt} \in \arg\min\limits_{P \in \Omega^N}J^N(P)$
and $P^* \in \arg\min\limits_{P \in \Omega}J(P)$.
For any $P \in \Omega^N$, the OT problem \eqref{eq: discrete Monge prime} 
has the same optimal value than its relaxation
\eqref{eq: discrete Kantorovich prime}, i.e., $\bar{C}^N(P)=C^N(P)$ \cite[Theorem 14.2]{Vanderbei2020LinearExtensions}.
Therefore, for all $P \in \Omega^N$, $J^N(P)=\bar{J}^N(P)$ and so in particular
\begin{equation} \begin{aligned}
\bar{J}^N \left( P_{opt}^N \right) = J^N \left( P_{opt}^N \right) = J^N_{opt}.
\nonumber \end{aligned}\end{equation}
Let $\varepsilon > 0$. By Lemma \ref{lemma: value function uniform convergence}, with probability one, 
$\bar{J}^N$ converges uniformly to $J$ as $N \to \infty$. 

Hence, there exists $N_{\varepsilon} \in \mathbb{N}$ such that for all $N \geq N_{\varepsilon}$, we have for all $P \in \Omega$
\[
J(P) \leq \bar J^N(P) + \varepsilon \text{ and } \bar J^N(P) \leq J(P) + \varepsilon.
\]
Hence we get, for all $N \geq N_\varepsilon,$
\begin{align*}
J_{opt} &= J(P^*) \leq J(P^N_{opt}) \leq \bar J^N(P^N_{opt}) + \varepsilon = J^N_{opt} + \varepsilon \\
\text{and } \; J^N_{opt} &= \bar J^N(P^N_{opt}) \leq \bar J^N(P^*) \leq J(P^*) + \varepsilon = J_{opt} + \varepsilon.
\end{align*}
This proves that with probability one, $\lim\limits_{N \rightarrow \infty} J^N_{opt} = J_{opt}$.

\subsection{Proof of Lemma \ref{lemma: conv sub-optimal cost to limit}}

\label{proof: conv sub-optimal cost to limit}

For any integer $N$ and any $P \in \Omega$, let $\tilde{\lambda}^N(P) \in \iset{D}^{N}$ be the vector of agent 
destination choices resulting from $\eqref{eq: mapping}$, i.e.,
$\tilde{\lambda}^N(P) \coloneqq \col(\tilde{\lambda}_1(P),\hdots,\tilde{\lambda}_N(P))$
with $\tilde{\lambda}_i(P) \coloneqq \lambda(x_i(0),P)$ for all $i \in \iset{N}$.
Note that $V^N(0,\tilde{\lambda}^N(P))$, with $V^N$ defined in \eqref{eq:simplified value function}, 
is equal to the social cost when the $N$ agents deploy the optimal control \eqref{eq:social optimal control} 
and use the agent-to-destination choices vector $\tilde{\lambda}^N(P)$. 

Let $P \in \Omega$.
To prove the lemma, we show that almost surely
$\tilde{J}^N(P)-V^N(0,\tilde{\lambda}^N(P))$ converges to $0$ and
$V^N(0,\tilde{\lambda}^N(P))$ converges to $J(P)$ as $N \to \infty$.

For the first part, and recall that $\tilde J^N$ is defined by \eqref{eq: tilde J}.
For any $i \in \iset{N}$, recalling the notation \eqref{eq: F^N},
let $x_{i}^{N,*}(t,F^N(P))$ 
be the trajectory of agent $i$ when it uses 

the optimal control $u_{\tilde{\lambda}_i(P)}^{N,*}$ from \eqref{eq:social optimal control}. 
We have 
\begin{align*}
&V^N(0,\tilde{\lambda}^N(P))= \\
&\frac{1}{N} \sum_{i=1}^N J_{i,\tilde{\lambda}_i(P)}\left(u^{N},*_{i,\tilde{\lambda}_i(P)}
\left(\cdot,x_{i}^{N,*}(\cdot,F^N(P)),F^N(P)\right), \right.\\
& \qquad \qquad \qquad \qquad \qquad \qquad \qquad \qquad \left. \bar{x}^N\left(\cdot,F^N(P)\right)\right).
\end{align*}

Agents are initially randomly distributed identically and independently according to the distribution $\mathcal{P}_0$. 
By the SLLN under Assumption \ref{asspt: E compact} and \eqref{eq: x_bar},  we get that almost surely
\begin{align}
\lim_{N \rightarrow \infty}\frac{1}{N}\sum_{i=1}^N x_i(0) = \int_{x \in E} x \, d \mathcal P_0 = \bar{x}(0,P). 
\label{eq: second result}
\end{align}
Moreover, again by the SLLN, $F^N(P)$ defined in \eqref{eq: F^N} 
converges almost surely to $P$ under Assumption \ref{asspt: measure zero}, 
because for the solution of the OT problem \eqref{eq:semi-discrete kantorovich prime}
we have $\mathcal P_0(\mathcal{C}_{j}^*(P)) = P_j$ for all $j \in \iset{D}$.
Using this result, \eqref{eq: second result}, and applying \cite[Theorem 3]{Coppel1965StabilityEquations} 
to the differential equations \eqref{eq: x_bar} and \eqref{eq: x_tilde}, we conclude that almost surely $t \to \tilde{x}^N(t,P)$ converges uniformly on $[0,T]$ to $t \to \bar{x}(t,P)$ as
$N \to \infty$.
Using similar arguments, almost surely the function $t \to \bar{x}^N(t,F^N(P))$, which satisfies the differential equation
\eqref{eq: mean dynamics} for $P^N = F^N(P)$, also converges uniformly on $[0,T]$ 
to $t \to \bar{x}(t,P)$.

Therefore, almost surely the function $t \to \bar{x}^N(t,F^N(P))-\tilde{x}^N(t,P)$ converges uniformly to $0$ on $[0,T]$.
Consider now the control laws \eqref{eq:social optimal control} and \eqref{eq:approx control}. 
By the same argument as in the proof of Lemma \ref{lemma: squared uniform convergence} above, 
the function $t \to \varphi_2^N(t) \bar{x}^N(t,F^N(P))$ converges almost surely uniformly on $[0,T]$ 
to  $t \to \varphi_2(t) \bar{x}(t,P)$. 
From this and Lemma \ref{lemma: alpha, beta, phi convergence}, 
for each $j \in \llbracket D \rrbracket$ and $x \in E$, the function
$t \to u_{j}^{N,*}(t,x,F^N(P))$ converges almost surely uniformly on $[0,T]$ to $t \to \tilde{u}_{j}(t,x,P)$
as $N \to \infty$.
Consequently, almost surely the closed-loop trajectories
$t \to x_{i}^{N,*}(t,F^N(P))$ converge uniformly on $[0,T]$ to $t \to \tilde{x}_i(t,P)$ as $N \to \infty$,
for $i \in \iset{N}$.

Let $\epsilon > 0$.
Let $j \in \llbracket D \rrbracket$. Combining the almost sure uniform convergence of agent state, 
control and population average trajectories, we conclude that for any agent of index $i$ 
whose final destination index is $j$, there exists $N_{j,}\varepsilon \geq 0$ such that for all 
$N \geq N_{j,}\varepsilon$, for all $t \in [0,T]$,
the costs in \eqref{eq:individual  cost} satisfy
\begin{align*}
&\left| \left[-[x_{i}^{N,*}(t,F^N(P))- \bar{x}^N(t,F^N(P))]_{R_x}^2 \right.\right.\\
&\left.\left.+[x_{i}^{N,*}(t,F^N(P)) - d_{j}]^2_{R_d}+[\tilde{u}_{j}(t,x_{i}^{N,*}(t,F^N(P)),P)]^2_{R_u}\right]\right.\\
&\left.- \left[-[\tilde{x}_i(t,P)- \tilde{x}^N(t,P)]_{R_x}^2\right.+[\tilde{x}_i(t,P) - d_{j}]^2_{R_d} \right.\\
& \left. \left.+[\tilde{u}_{j}(t,\tilde{x}_i(t,P),P^N)]^2_{R_u}\right] \right| \leq \varepsilon,
\end{align*}

and
\[
\left| \left[ x_{i}^{N,*}(T,F^N(P)) - d_j \right]^2_M - \left[ \tilde x_{i}(T,P) - d_j \right]^2_M \right|
\leq \varepsilon.
\]
Take $N_{\varepsilon}=\max_{j \in \llbracket D \rrbracket} N_{j},\varepsilon$. 
This shows that for all $i \in\iset{N}$
\begin{align*}
& \Big |J_{i,\tilde{\lambda}_i(P)}(\tilde{u}_{\tilde{\lambda}_i(P)}\left(\cdot,\cdot,P\right),\tilde{x}^N\left(\cdot,P\right))\\
&-J_{i,\tilde{\lambda}_i(P)} \left( u^{N,}*_{i,\tilde{\lambda}_i(P)}\left(\cdot,\cdot,F^N(P)\right) ,\bar{x}^N\left(\cdot,F^N(P)\right)\right) \Big|  \\
& \leq \varepsilon(T +1 ) .
\end{align*}
Taking the average, 
we get that 
almost surely

$\tilde{J}^N(P)-V^N(0,\tilde{\lambda}^N(P)) \to 0$ as $N \to \infty$,
which concludes the first step of the proof.
Next, we prove that $V^N(0,\tilde{\lambda}^N(P))$ converges almost surely to $J(P)$, defined by
\eqref{eq:simplified value function infinite}. Using \eqref{eq:simplified value function}, 
we have 
\begin{align}
&V^N(0,\tilde{\lambda}^N(P))=\frac{1}{2N}\sum\limits_{i=1}^{N}x_i^T(0)\left(\varphi^N_1(0)-\frac{\varphi^N_2(0)}{N}\right)x_i(0) \nonumber \\
&+\frac{1}{2}(\bar{x}^N(0))^T\varphi^N_2(0)\bar{x}^N(0) + 
\frac{1}{2N}\sum\limits_{i=1}^{N}|x_{i}(0)-\beta^N_{\tilde{\lambda}_i(P)}(0)|_2^2 \nonumber \\
&+\sum\limits_{j=1}^{D-1}F^N_j(P)(\alpha_j^N(0))^T\bar{x}^N(0)-\frac{1}{2}\sum\limits_{j=1}^{D}F^N_j(P)|\beta_{j}^N(0)|_2^2 
\nonumber \\
& -\frac{1}{2N}\sum\limits_{i=1}^{N}|x_i(0)|_2^2 + \chi^N(0,F^N(P)). \label{eq: tildeV^N(0)}
\end{align}
By the definition of $\lambda(x,P)$ in \eqref{eq: mapping}, we can rewrite the term
\begin{align*}
&\frac{1}{N} \sum_{i=1}^{N} \left| x_i(0)-\beta^N_{\tilde{\lambda}_i(P)}(0) \right|_2^2 = \\
&= \frac{1}{N} \sum_{i=1}^{N} \sum_{j=1}^D \left| x_i(0)-\beta^N_{j}(0) \right|_2^2 \; 1_{\{ x_i(0) \in \mathcal{C}_{j}^*(P)\}}
\\
&= \frac{1}{N}\sum_{i=1}^{N} | x_i(0)|_2^2 +\sum_{j=1}^{D}F_j^N(P) \left|\beta^N_{j}(0) \right|_2^2\\
& \quad -2\sum_{j=1}^{D} \beta^N_{j}(0)^T\left(\frac{1}{N}\sum\limits_{i=1}^N x_i(0) 
\; 1_{ \{ x_i(0) \in \mathcal{C}_{j}^*(P)\}} \right).
\end{align*}
Using the SLLN under Assumption \ref{asspt: E compact},
the fact that $F^N(P)$ converges almost surely to $P$ under Assumption \ref{asspt: measure zero}
as explained below \eqref{eq: second result}, and Lemma \ref{lemma: SDOT result}, we conclude that almost surely,
\begin{align*}
&\lim_{N \rightarrow \infty} \frac{1}{N}\sum_{i=1}^{N}
\left| x_i(0)-\beta^N_{\tilde{\lambda}_i(P)}(0) \right|_2^2  
=\int_{x \in E} |x|_2^2 \, d \mathcal{P}_0 \\
& \quad +\sum_{i=1}^{D}P_j|\beta_{j}(0)|_2^2
-2\sum_{j=1}^{D}\beta_{j}(0)^T\int_{x \in E} x \, 1_{ \{ x \in \mathcal{C}_{j}^*(P) \} } d \mathcal{P}_0\\
&=\int_{x \in E} \left| x-\beta_{\lambda(x,P)}(0) \right|_2^2  d \mathcal{P}_0 = C(P).
\end{align*}

Then, using again Lemma \ref{lemma: alpha, beta, phi convergence}, the SLLN, and the fact
that $F^N(P)$ converges almost surely to $P$, 
it is straightforward to show the almost sure convergence of
the other terms in \eqref{eq: tildeV^N(0)} to the corresponding terms
in \eqref{eq:simplified value function infinite} as $N \to \infty$, which concludes the proof.
\subsection{ Proof of Lemma \ref{Lemma: subgradient}}
\label{proof: subgradient}
Using the expressions \eqref{eq:simplified value function infinite} and \eqref{eq: chi polynomial infinite}, we show that,
\begin{equation}
\begin{aligned}
\label{eq: partial gradient}
   &\nabla_P \left(J(P)-\frac{1}{2}C(P)\right)=2\int_T^0 WdtP\\
   &+\frac{1}{2}\intop_T^0\begin{pmatrix} 
      [\beta_1]^2_S -[d_1]^2_{R_d}\\
      \vdots\\ 
      [\beta_{D-1}]^2_S -[d_{D-1}]^2_{R_d}\\
       [\beta_D]^2_S -[d_D]^2_{R_d}
   \end{pmatrix}dt\\
   &+\begin{pmatrix} 
      \alpha_1(0)^T\bar{x}(0)+\frac{1}{2}([d_1]^2_{M}-|\beta_{1}(0)|_2^2)\\
      \vdots\\ 
      \alpha_{D-1}(0)^T\bar{x}(0)+\frac{1}{2}([d_{D-1}]^2_{M}-|\beta_{D-1}(0)|_2^2)\\
      \frac{1}{2}([d_D]^2_{M}-|\beta_{D}(0)|_2^2)
   \end{pmatrix}.
\end{aligned}
\end{equation}
Moreover, for all $P_1, P_2 \in \Omega$, we have,
\begin{equation}
\begin{aligned}
& C(P_2)=J^D(g(P_2),P_2) \geq J^D(g(P_1),P_2).
\end{aligned}
\end{equation}
Hence,
\begin{equation}
\begin{aligned}
&C(P_2) \geq \sum_{j=1}^D \int\limits_{x \in E} \min\limits_{j \in 1,\hdots,D}(|x-\beta_j(0)|_2^2-g_j(P_1)) d \mathcal{P}_0(x)\\
&\qquad \qquad \qquad +\left\langle g(P_1), P_1\right\rangle+\left\langle g(P_1), P_2-P_1\right\rangle \\
&\qquad \qquad \geq C\left(P_1\right)+\left\langle g(P_1), P_2-P_1\right\rangle,
\end{aligned}
\end{equation}
with $\left\langle.,.\right\rangle$ is the dot product.
Consequently, for all $ P \in \Omega, \; g(P)$ is a subgradient of $C(P)$ at $P$.

 Therefore, using the latter result and \eqref{eq: partial gradient}, we get the result.

\end{document}